\documentclass[12pt]{article}
\usepackage{srcltx}
\usepackage{amssymb}
\usepackage{amsfonts}
\usepackage{latexsym}
\usepackage{amsmath}
\usepackage{amsthm}
\usepackage{graphicx}

\newcommand{\refLem}[2]{ Lemma~\ref{#2}.\ref{#1} }
\renewcommand{\theequation}{\arabic{section}.\arabic{equation}}

\theoremstyle{plain}
\newtheorem{theorem}{\bf Theorem}[section]
\newtheorem{lemma}[theorem]{\bf Lemma}
\newtheorem{corollary}[theorem]{\bf Corollary}

\newtheorem{definition}[theorem]{\bf Definition}
\theoremstyle{remark}

\def\[{\begin{equation}}
\def\]{\end{equation}}
\let\ge\geqslant
\let\le\leqslant

\def\Re{\mathop{\rm Re}\nolimits}
\def\Im{\mathop{\rm Im}\nolimits}

\def\sign{\mathop{\rm sign}\nolimits}
\def\l{\lambda}
\def\ve{\varepsilon}
\newcommand{\co}{C}
\newcommand{\Ai}{\hbox{\rm Ai}}
\newcommand{\Bi}{\hbox{\rm Bi}}
\newcommand{\ai}{ a}

\newcommand{\re}{\Re}
\newcommand{\vt}{\vartheta}
\newcommand{\vk}{\varkappa}
\newcommand{\C}{\mathbb{C}}
\newcommand{\ol}{\overline}

\def\zlx{z_\l(x)}
\def\zlpx{z_\l'(x)}

\begin{document}
\markboth{Estimates of parabolic cylinder functions
}{Estimates of parabolic cylinder functions}

\title{Estimates of parabolic cylinder functions\\
on the real axis}
\author{
  A.A. POKROVSKI\\
    Laboratory of Quantum Networks, Institute for Physics,\\
    St-Petersburg State University, St.Petersburg~198504,\\
    Ulyanovskaya~1.\\
    E-mail: alexis.pokrovski@mail.ru
}
\date{}
\maketitle

\renewcommand{\abstractname}{}
\begin{abstract}
We estimate the  expressions $F(\pm x,\l)$ and $F(\pm ix,-\l)$, where
$F(x,\l)=\frac{\partial}{\partial x}U(-\tfrac{\l}{2}, x\sqrt{2})+U(-\tfrac{\l}{2}, x\sqrt{2})\sqrt{x^2-\l}$
and  $U$ is the
standard solution of the parabolic cylinder equation, satisfying
$U(a,x)\sim x^{-a-1/2} e^{-x^2}$ as $x\to+\infty$.
The estimates are valid in rather complicated domains and
refine there the classical result of Olver.
The estimates with real $x$ are important for the
 spectral analysis of non-analytically perturbed quantum harmonic oscillator.
We determine the part of the real  $x$-axis, which is  within the domains of the estimates.
This requires a detailed study of the image of the real $x$-axis in the
standard quasiclassical variable.


\end{abstract}

\section{Introduction}

Consider the solution $\psi$
of
the  equation
\begin{equation}
 \label{OurEq}
-y''(x,\l)+x^2y(x,\l)=\l y(x,\l),
\qquad
y'\equiv\frac{\partial y}{\partial x},
\end{equation}
defined by
  $\psi(x,\l)=U(-\tfrac{\l}{2}, x\sqrt{2})$,
  where $U$ is the standard parabolic cylinder function
 $U(a,x)=D_{-a-\frac{1}{2}}(x)$
(see
\cite{AS},\cite{WW}).
Our main result is the estimates of the  expressions
$F(\pm x,\l)$ and  $F(\pm ix,-\l)$, where
$F( x,\l)=\psi'( x,\l)+\psi( x,\l)\sqrt{x^2-\l}$.
Each of the four $F$'s is evidently related to one of the four  solutions
$\psi(\pm x,\l)$,  $\psi(\pm i x,-\l)$ of (\ref{OurEq}).
The estimates and
their (rather complicated) domains are
given in Theorem~\ref{Basic:InVariableX}.
In particular, for $F( x,\l)$ the estimate has the form
\begin{equation}\label{Example}
|F(x,\l)|
\le
\co
|\phi(\l)
\rho(x,\l)
|e^{-\l\xi(\frac{x}{\sqrt{\l}})}|
,
\qquad
\rho(x,\l)=
\frac{1+|\l|^{\frac{1}{2}}+|x^2-\l|^{\frac{1}{2}}}{1+|\l|^{\frac{5}{12}}+|x^2-\l|^{\frac{5}{4}}}
\end{equation}
where $\co$ is an absolute constant,
 $\xi(t)=\int\limits_1^t\!\!\!\sqrt{s^2-1}\,ds$ and
$\phi(\l)=
2^{\frac{3}{4}}\sqrt{\pi}\left(\frac{\l}{2e}\right)^{\frac{\l}{4}}$ are
positive on $(1,\infty)$ and defined on
$\mathbb{C}\setminus(-\infty,1]$
 and $\mathbb{C}\setminus \mathbb{R}_-$, respectively.
The estimates  of $F(-x,\l)$ and  $F(\pm ix,-\l)$ have similar form,
with modifications taking into account branching
of $\xi$, $\phi$ and $\sqrt{x^2-\l}$.

We take special care of the case of real $x$; to this end
we study the image of the real axis in the relevant quasiclassical variable
(see Section~\ref{Section:PropGammaL}).
Then we localize the image within the domains of the estimates.
By Theorem~\ref{Basic:InVariableX}, for each real $x$ and each ray $\arg \l=$const the
estimates for at least two of the four expressions
$F(\pm x,\l)$,  $F(\pm ix,-\l)$ are fulfilled and the two related
 solutions  $\psi$ are linearly independent. The last issue is important for applications that we discuss below.

Our result is a refinement
of the following estimates  due to Olver
\cite{Olver1960}:
\begin{equation}\label{OlversEstimates}
|\psi(x,\l)|
\le
\co
\frac{|\phi(\l)|}{\rho_0(x,\l)}
|e^{-\l\xi(\frac{x}{\sqrt{\l}})}|,
\qquad
|\psi'(x,\l)|
\le
\co
|\phi(\l)
\rho_0(x,\l)
e^{-\l\xi(\frac{x}{\sqrt{\l}})}|,
\end{equation}
where $\rho_0(x,\l)=1+|\l|^{\frac{1}{12}}+|x^2-\l|^{\frac{1}{4}}$.
These estimates are valid for
all values of $x$ and $\l$ (for special choice of branches of
$\xi$, $\phi$ and $\sqrt{x^2-\l}$, discussed in \cite{Olver1960}),
whereas the domain of our estimates is smaller. For example,
for  $\l>0$ the asymptotics (\ref{Example}) holds true only for
$\sqrt{\l}\le x$.

The results may have little interest in their own right,
but they are useful for  the study of
the quantum harmonic oscillator, perturbed by a non-decaying potential
(\hbox{$-\frac{d^2}{dx^2}+x^2+q(x)$} on $L^2(\mathbb{R})$,
where $q(x)$ is not small for large $|x|$).
For perturbations $q$ satisfying
$\sup_{x\in\mathbb{R}}(|q'(x)|+|\int_0^xq(s)\,ds|)<\infty$
the following  spectral asymptotics
was proved in \cite{ArXiv0}
\begin{equation}\label{SpectralAsymptotics}
\mu_n=
\mu_n^0+
\frac{1}{2\pi}
\int_{-\pi}^\pi
q(\sqrt{\mu_n^0}\sin\vartheta)
d\vartheta
+
O(n^{-\frac{1}{3}}),
\qquad
\mu_n^0=2n+1.
\end{equation}
The estimates, obtained in the present paper, give a tool for improvement
of the error bound in (\ref{SpectralAsymptotics}).
We plan to make this in a separate paper.

%

In the proof we follow the classical scheme used in \cite{Olver1959}.
In Section~\ref{Section:ElemProrParabCylFunc}
we list the necessary properties of the parabolic cylinder functions.
In Section~\ref{Section:ChangeOfVar} we introduce
 the  quasiclassical variable
$z=z_\l(x) $; the
equation (\ref{OurEq}) becomes a
perturbed Airy equation (\ref{EqinZ}) with
the perturbation $V_0$ 
decaying in both
$z$ and $\l$.
In Section~\ref{Section:PropGammaL} we study the properties of the
family of curves
$\Gamma_\l=z_\l(\mathbb{R}_+)$, taking special care of its positioning
relatively to the sectors of decay of  Airy functions. Since $\Gamma_\l$
for small $\arg\l$ is not within one such sector,
we split the curve by the image
  $z_*=z_\l(x_*)$ of a suitably defined turning point $x_*\in \mathbb{R}_+$.

In Section~\ref{Section:FinalEstimates}
we formulate the main result. We prove it first in $z$-variable.
We fix four solutions $A_0$, $A_\pm$ and $A_*$ of Eq.(\ref{EqinZ}),
 asymptotically close to  one of
the Airy functions.  For each solution we write
 integral equation.  Then we modify these integral equations by separating
the exponential multiplier, writing
$
A_0(z,\l)=e^{-\frac{2}{3}z^{\frac{3}{2}}} a_0(z,\l)
$
and similarly introducing $a_{\pm}$, $a_*$.
Analyzing the modified integral equations,
 we prove the main estimates in terms
of $z$-derivatives of
$a_0$, $a_\pm$ and $a_*$ (Theorem~\ref{Basic:Svoistva}).
The domains of these estimates turn out to be rather complicated.
Comparing the asymptotics as $x\to\infty$, we identify each of
 $A_\nu(z_\l(x),\l)/\sqrt{z_\l'(x)}$
with one of $\psi(\pm x,\l)$,
$\psi(\pm ix,-\l)$.
This yields
Theorem~\ref{Basic:InVariableX}. We also give the connection
formulas  and calculate the Wronskians  for $A_\nu$.

In Appendix~A we list the properties of the auxiliary family of curves,
 that are used in the proof of Theorem~\ref{Basic:Svoistva}.
In Appendix~B we accomplish the study of properties of
$\Gamma_\l=z_\l(\mathbb{R}_+)$
by estimating the integrals of
$\frac{1}{(1+|z|)^\alpha}$
and
$\frac{|e^{\pm z^{\frac{3}{2}}}|}{(1+|z|)^{\alpha}}$ along the family of
curves $\Gamma_\l(z)$.
For the sake of completeness  some
technical
results from \cite{ArXiv0} are reproduced
in Section~\ref{Section:PropGammaL}
and Apendix~B.

\textbf{Notations.}
\begin{itemize}
    \item We set
$\mathbb{C}_+=\{z\in\mathbb{C}:\Im z>0\}$
,
$\overline{\mathbb{C}}_+=\{z\in\mathbb{C}:\Im z\ge 0\}$.
    \item For the spectral parameter $\l\in \overline{\C}_+\setminus \{0\}$ we set
    $\l=|\l|e^{2i\vt}$,
    $\vt\in[0,\frac{\pi}{2}]$.
    \item The functions $\log z$  and $z^\alpha=e^{\alpha\log z}$
    for $\alpha\in\mathbb{C}$
take their principal values on
$\mathbb{C}\setminus\mathbb{R}_-$.
    \item We denote by $S[\alpha,\beta]$  the sector
$\{z\in\mathbb{C}: \arg z\in[\alpha,\beta]\}$; similarly we define
$S(\alpha,\beta)$, $S[\alpha,\beta)$ and $S(\alpha,\beta]$.
By $S[-\pi,\pi]$ we denote the complex plane  cut along $(-\infty,0]$,
where the upper and the lower sides of the cut are included,
but not identified.
    \item We also set
$S_R[\alpha,\beta]=\{z\in S[\alpha,\beta]:|z|\ge R \}$.
\end{itemize}

\section{Elementary properties of parabolic cylinder functions}
\label{Section:ElemProrParabCylFunc}

In this section we list  the necessary properties of
standard parabolic cylinder functions. We use
the solution $\psi$ of Eq.(\ref{OurEq}), given by
 $\psi(x,\l)=U(-\tfrac{\l}{2}, x\sqrt{2})$, where $U$
is the standard parabolic cylinder function (see
\cite{AS}). In the notation of Whittaker \cite{WW}
for the parabolic cylinder function $D_n(z)$
and the confluent hypergeometric function $W_{k,m}(x)$
 we have
$
\psi(x,\l)=D_{\frac{\l-1}{2}}(x\sqrt{2})=
\frac{2^{\frac{\l-1}{4}}}{\sqrt{x}}W_{\frac{\l}{4},-\frac{1}{4}}(x^2).
$
The solution $\psi$
 is uniquely defined by its asymptotics
\begin{equation}
\label{DefinitionParabolic CylFunct}
\psi(x,\l)\sim (x\sqrt{2})^{\frac{\l-1}{2}}e^{-\frac{x^2}{2}}
\quad
\hbox{as} \quad
|x|\to\infty,
\quad
 |\arg x|\le\tfrac{3\pi}{4}- \epsilon \quad\hbox{for any} \quad\epsilon>0.
\end{equation}
It is an entire function of $x$ and an entire function of $\l$.
Other important solutions of (\ref{OurEq}) are
$\psi(- x,\l)$ and $\psi(\pm ix,-\l)$. The connection formulas are
\begin{align}\label{ConnFormIX}
\psi(\pm ix,-\l)
&=
\frac{\Gamma(\frac{1-\l}{2})}{\sqrt{2\pi}}
\left(
e^{i\frac{\pi}{4}(\l+1)}\psi(\pm x,\l)+
e^{-i\frac{\pi}{4}(\l+1)}\psi(\mp x,\l)
\right)
\\
\label{ConnFormMinusX}
\psi(\pm x,\l)
&=
\frac{\Gamma(\frac{1+\l}{2})}{\sqrt{2\pi}}
\left(
e^{i\frac{\pi}{4}(\l-1)}\psi(\pm ix,-\l)+
e^{-i\frac{\pi}{4}(\l-1)}\psi(\mp ix,-\l)
\right)
\end{align}

\section{The changes of variables}
\label{Section:ChangeOfVar}

In the variable $t=\frac{x}{\sqrt{\l}}$ the equation (\ref{OurEq})
becomes
\begin{equation}\label{EqTCoord}
w''(t)=\l^2(t^2-1)w(t),
\qquad
w(t)=y(t\sqrt{\l} ).
\end{equation}
Introduce the function
$\xi(t)=\int\limits_1^t\!\!\!\sqrt{s^2-1}\,ds$, such that $\xi>0$ for $t>0$,
defined on
$\mathbb{C}\setminus(-\infty,1]$.
We have
\begin{equation}
\label{DefXi}
\xi(t)
=
\tfrac{1}{2}\left(t\sqrt{t^2\!-1}-\log(t+\sqrt{t^2-1})\right)
=
\tfrac{t^2}{2}-\tfrac{1}{2}\ln 2t- \tfrac{1}{4}+O(t^{-2}),
\qquad t\to\infty.
\end{equation}
$\xi$
has finite branch points $t=+1$ and $t=-1$,
which we denote by $A$ and $E$, respectively.

{\bf Definition of $\eta(t)$. }
Introduce the function
\begin{equation}\label{DefEta}
\eta(t)=\left(\tfrac{3}{2}\xi(t)\right)^{\frac{2}{3}},
\quad
t\in \C\setminus(-\infty,-1].
\end{equation}
Unlike $\xi$, it is analytic on $(-1,1)$ and $t=1$
is its regular point.
 Define the points
$\eta_0=\eta(0)=-(\frac{3\pi}{8})^{\frac{2}{3}}$
and
$\eta_E=\eta(-1)=-(\frac{3\pi}{4})^{\frac{2}{3}}$.
 The mapping
$
\eta: D_T\mapsto\C\setminus(-\infty,\eta_E]
$ is an analytic isomorphism, where
$
D_T=(-1,1]\cup
\{
t: |\arg \xi(t)|<\tfrac{3\pi}{2}
\}
$.
  We denote by $t(\eta)$ the function inverse to
$\eta(t)$ and defined on $\C\setminus(-\infty,\eta_E]$.

\begin{figure}[ht]
\centering 
\includegraphics[width=165mm]{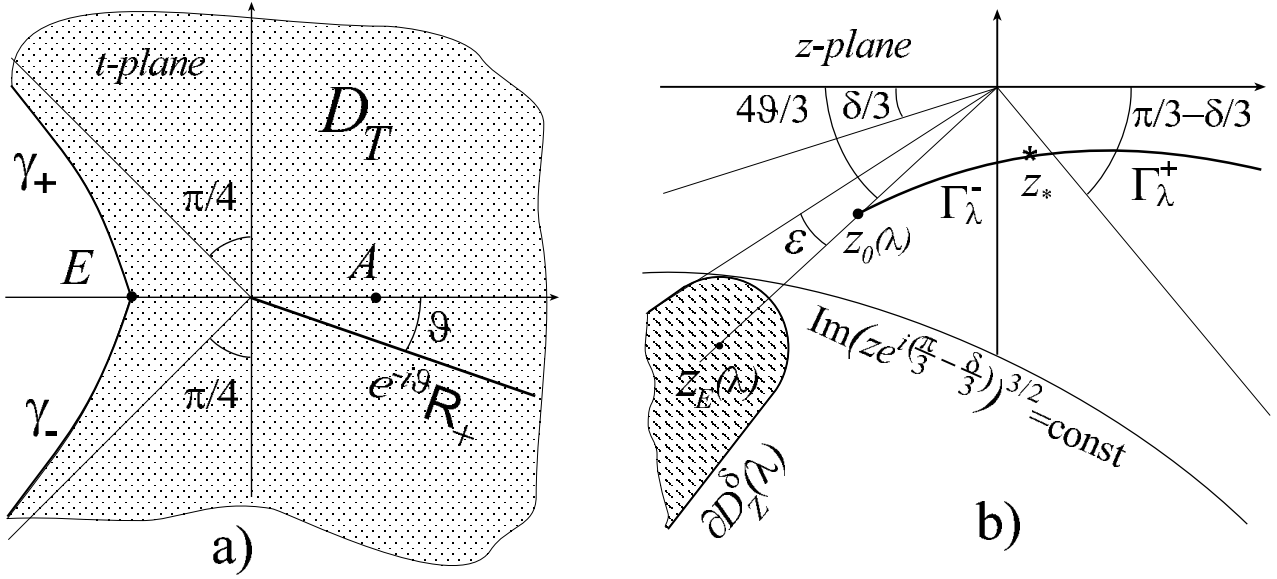}\\
\caption[short caption here]{
a) The ray $e^{-i\vt}\mathbb{R}_+$ for
 $\vt\in [0,\frac{\pi}{2}]$ is in the domain $D_T$ in $t$-plane.
b) The point $z_*$ divides the curve
$\Gamma_\l=\Gamma_\l^-\cup\Gamma_\l^+$.  The curve
$\Im (ze^{i(\frac{\pi}{3}-\frac{\delta}{3})})^{\frac{3}{2}}=${\rm const} is
tangential to the boundary of $D_Z^\delta(\l)$ and is separated away
from $\Gamma_\l$. The complement to $D_Z^\delta(\l)$
is dashed. Here $2\vt=\arg\l$.
}
\label{DomainAndGamma+DZ}  
\end{figure}

The boundary of $D_T$ is $\gamma_-\cup\gamma_+$,
where
$
\gamma_\pm=\{t\in\overline{\C}_\pm:\arg\xi(t)=\pm\frac{3\pi}{2},
|\xi(t)|\ge\frac{\pi}{2}\}
$.
The curves $\gamma_\pm$ and $\gamma_{E+}$ emanate from  $E$
and have tangential inclinations to the positive real axis
of $\pm\frac{2\pi}{3}$, respectively.
By (\ref{DefXi}), the curves $\gamma_\pm$  are asymptotic to
the rays
$
\arg t=\pm\frac{3\pi}{4}
$,
respectively.
The domain $D_T$ and the curves $\gamma_\pm$ are
 schematically presented on Fig.\ref{DomainAndGamma+DZ}~a).

For any fixed
$\l\in S[-\pi,\pi]\setminus\{0\}$
transition from $x$-variable to $z$-variable is given by
 the function
\begin{equation}\label{DefZofX}
z_\l(x)=\l^{\frac{2}{3}}\eta(\tfrac{x}{\sqrt{\l}})
\equiv
\l^{\frac{2}{3}}
\left(\tfrac{3}{2}\xi(\tfrac{x}{\sqrt{\l}})\right)^{\frac{2}{3}},
\quad
\quad
x\in D_X(\l)
\overset{\text{def}}{=}
\{x\in\C: \tfrac{x}{\sqrt{\l}}\in D_T\}.
\end{equation}
For a fixed $\l\in S[-\pi,\pi]\setminus\{0\}$ each mapping
\begin{equation}\label{DefDomainDZ}
z_\l(\cdot):D_X(\l)\mapsto D_Z(\l)\overset{\text{def}}{=}
\C\setminus\{z:\arg z=\arg z_E(\l),
|z|\ge|z_E(\l)|\}
\end{equation}
is an analytic isomorphism. We denote by $x_\l(\cdot)$ the function inverse to
$z_\l(\cdot)$ and defined on $D_Z(\l)$.
Evidently $z_\l(\sqrt{\l})=0$;
the images of $x=0$ and $x=-\sqrt{\l}$ in $z$-variable are
\begin{equation}
\label{Defz0zE}
z_0=z_\l(0)=\l^{\frac{2}{3}}\eta_0=
-\l^{\frac{2}{3}}\left(\tfrac{3\pi}{8}\right)^{\frac{2}{3}},
\quad
z_E(\l)=z_\l(-\sqrt{\l})=\l^{\frac{2}{3}}\eta_E=
-\l^{\frac{2}{3}}\left(\tfrac{3\pi}{4}\right)^{\frac{2}{3}}.
\end{equation}

{\bf The parabolic cylinder equation in $z$-variable.}
$y$ is a solution of (\ref{OurEq}) in $D_X(\l)$ if and only if
$u(z,\l)=\frac{y(x_\l(z),\l)}{\sqrt{\partial_z x_\l(z)}}$
solves the equation
\begin{equation}
\partial_z^2 u(z,\l)-zu(z,\l)= V_0(z,\l)u(z,\l),
 \qquad z\in D_Z(\l),
  \label{EqinZ}
\end{equation}
where we use the notation $\partial_z u=\frac{\partial u}{\partial z}$ and where
  \begin{equation}
  \label{EffPotentials}
  V_0(z,\l)={v\left({z}{\l^{-\frac{2}{3}}}\right)}{\l^{-\frac{4}{3}}},
  \quad
v(\eta)=\sqrt{t'(\eta)}\frac{d^2}{d\eta^2}\frac{1}{\sqrt{t'(\eta)}},
\end{equation}
and $t(\cdot)$ denotes the function, inverse to $\eta(\cdot)$,
given by (\ref{DefEta}).
The function $v(\eta)$ is analytic in $\C\setminus(-\infty,\eta_E]$,
hence
 $V_0$ is analytic in both $\l$ and $z$ for
 $(\l,z)\in(\mathbb{C}\setminus\mathbb{R}_-)\times D_Z(\l)$.
Since $v(\eta)$ has the uniform asymptotics
$
v(\eta)\sim\frac{7}{64}\eta^{-2}
$ as $|\eta|\to\infty$, it is unbounded only at $\eta=\eta_E$.
 Thus for any $\epsilon>0$  there exists a constant $C$ such that
\begin{equation}
 \label{Estv0-eta}
  |v(\eta)|\le\frac{\co}{(1+|\eta|)^2},
  \qquad
  \text{for}
  \quad
  |\eta-\eta_E|\ge\epsilon,
  \quad
\eta\not\in
(-\infty,\eta_E]
.
\end{equation}

{\bf The curve $\Gamma_\l=z_\l(\mathbb{R}_+)$. }
Each mapping $z_\l(\cdot):\mathbb{R}_+\to
\Gamma_\l=z_\l(\mathbb{R}_+)$ is a real analytic isomorphism.
 If $\l>0$, then
$\Gamma_\l=[z_0,\infty)\subset\mathbb{R}$ is a half-line ($z_0$ is
given by (\ref{Defz0zE})). The curve is schematically presented on
  Fig.~\ref{DomainAndGamma+DZ} b).
For $z_1,z_2\in\Gamma_\l$ such that $x_\l(z_1)\le x_\l(z_2)$
 define the curves,
 that play the role of an interval, by
$$
\Gamma_\l(z_1,z_2)=\{z: x_\l(z)\in[x_\l(z_1),x_\l(z_2)]\},\ \ \ \
\Gamma_\l(z_1)\equiv \Gamma_\l(z_1,\infty)=\{z: x_\l(z)\ge x_\l(z_1)\}.
$$

Now we generalize the notion of the turning
point for non-positive $\l$.
\begin{definition}
\label{Z-star}
For each $\l\in S[-\pi,\pi]\setminus \{0\}$
define $x_*\equiv x_*(\l)\in\mathbb{R}_+$ and
$z_*\equiv z_*(\l)\in\Gamma_\l$ by
$|z_\l(x_*)|=\min\limits_{x\in \mathbb{R}_+}\ |z_\l(x)|$
and $z_*=z_\l(x_*)$.
We also define  $t_*$ and $r_*$ by
  $t_*=\frac{x_*}{\sqrt{\l}}$ and $r_*=|t_*|$.
\end{definition}
By Lemma~\ref{DivG2}, the point $z_*$ is defined correctly (is
unique). Throughout the paper  we use $z_*$ according to this
definition, omitting dependence from $\l$ for brevity.  We set
\begin{equation}
\label{DefGPM}
\Gamma_\l^-=\Gamma_\l(z_0,z_*),\qquad
\Gamma_\l^+=\Gamma_\l(z_*,\infty),\qquad
\Gamma_\l=\Gamma_\l^-\cup\Gamma_\l^+.
\end{equation}
Note that for $\l>0$ we have $x_*=\sqrt{\l}$.



\textbf{Separating small and large $\arg \l$. }
Further analysis of Eq.(\ref{EqinZ}) employs different technique
for small and large arguments of $\l$.  We formalize these
two cases by introducing  $\delta$ and considering separately
$|\arg\l|\le\delta$ and $\delta<|\arg\l|\le\pi$.
Here and below we fix
\begin{equation}\label{DefDelta}
\delta\in(0,\frac{\pi}{5}).
\end{equation}
Throughout the paper the constant $\co$ in the estimates depends on $\delta$;
for brevity we indicate this dependence explicitly only in
Theorem~\ref{Basic:InVariableX}.

For each $\l$ define the  domain
$D_Z^\delta(\l)$
in $z$-plane as follows.
Let $B_\ve(\l)$ denote the
disk of radius
$|z_E(\l)|\sin\ve$, centered at $z_E(\l)$.
Define the points $w_{\pm\delta}(\l)$ where
the curves
$\Im (ze^{i(\frac{\pi}{3}\mp\frac{\delta}{3})})^{\frac{3}{2}}=$const
are tangential to $D_\ve(\l)$:
 $\Im (w_{\pm\delta}(\l)
e^{i(\frac{\pi}{3}\mp\frac{\delta}{3})})^{\frac{3}{2}}
=
\sup_{z\in B_\ve(\l)}
\Im (ze^{i(\frac{\pi}{3}\mp\frac{\delta}{3})})^{\frac{3}{2}}$.
$w_\delta(\l)$ is defined for $\l\in S_{1/2}[\delta,\pi]$,
$w_{-\delta}(\l)$ is defined for $\l\in S_{1/2}[0,\pi-\delta]$.
 Fix $\ve\in(0,\frac{\delta}{3})$  sufficiently small
to ensure that for any
$\l\in S_{1/2}[\delta,\pi]$
 holds
$\Im (w_\delta(\l)e^{i(\frac{\pi}{3}-\frac{\delta}{3})})^{\frac{3}{2}}\le
\Im (z_0(\l)e^{i(\frac{\pi}{3}-\frac{\delta}{3})})^{\frac{3}{2}}$
\footnote{This choice of $\ve$ is used in Section~\ref{Section:FinalEstimates}
when we show that $\Gamma_\l$ (or its parts
$\Gamma_\l^{\pm}$) is within the range of
the estimates of Theorem~\ref{Basic:Svoistva}.}.
The domain's complement is the disk $B_\ve(\l)$ plus the disk's shadow from the point
light source at the origin,
\begin{equation}\label{DefD_Z^Epsilon}
D_Z^\delta(\l)=
\mathbb{C}\setminus
\left(
B_\ve(\l)
\cup
\{z:|\arg z-\arg z_E(\l)|\le \ve, |z|\ge|z_E(\l)|\cos\ve\}
\right),
\end{equation}

The complement to  $D_Z^\delta(\l)$ is schematically presented
on Fig.\ref{DomainAndGamma+DZ}~b) by the dashed region.

Here is the motivation of this definition and of the following one.
We prove the main estimate in $z$-variable using
integral equation, equivalent to (\ref{EqinZ}).
 The kernel of the integral equation
is a product of Airy functions and of the effective potential $V_0$.
By (\ref{Estv0-eta}), we have the estimate
\begin{equation}\label{EstV0-z}
 |V_0(z,\l)|\le \frac{C}{|\l|^{\frac{4}{3}}+|z|^2}
  \qquad
  \text{for}
  \quad
  z\in D_Z^\delta(\l).
\end{equation}
Airy functions allow
convenient estimates on the family of curves
$\Im (ze^{i\varphi})^{\frac{3}{2}}=${\rm const},
$|\varphi|\le \frac{\pi}{3}-\frac{\delta}{3}$
(the curves are discussed in details in Appendix~A).
Thus the estimates for the whole kernel hold in
a subdomain of $D_Z^\delta(\l)$,
whose points are attainable along a curve from this family.
Next we define an important part of this subdomain.

\begin{definition}
\label{DefHDelta}
$H_{\pm\delta}(\l)$ is the part of the sector
 $S[-\pi\pm\frac{\delta}{3},-\frac{\pi}{3}\pm\frac{\delta}{3} ]$,
 that lies above the curve
$\Im (ze^{i(\frac{\pi}{3}\mp\frac{\delta}{3})})^{\frac{3}{2}}=${\rm const},
tangential to the boundary of $D_Z^\delta(\l)$,
$$
H_{\pm\delta}(\l)=
\{
z\in S[-\pi\pm\tfrac{\delta}{3},-\tfrac{\pi}{3}\pm\tfrac{\delta}{3} ]:
\Im (ze^{i(\frac{\pi}{3}\mp\frac{\delta}{3})})^{\frac{3}{2}}
\ge
\Im (w_{\pm\delta}(\l)e^{i(\frac{\pi}{3}\mp\frac{\delta}{3})})^{\frac{3}{2}}
\}.
$$
$H_{+\delta}(\l)$ is defined for $\l\in S_{1/2}[\delta,\pi]$,
$H_{-\delta}(\l)$ is defined for $\l\in S_{1/2}[0,\pi-\delta]$.
\end{definition}

\section{Properties of $\Gamma_\l$}
\label{Section:PropGammaL}
In this section we find out the properties of the curve
$\Gamma_\l=z_\l(\mathbb{R}_+)$.
Since  $\xi$ is symmetric with respect to the real axis,
we consider  only  the case
$\Im\l\ge 0$. In this case we set $\xi(t)=\xi(t-i0)$ for
$t\in(-\infty,1]$ and follow this agreement throughout the paper.

 The relation
\begin{equation}\label{Z3-2PowInXi}
  \tfrac{2}{3}z_\l(x)^{\frac{3}{2}}=\l\xi (t),
  \qquad
  t=\tfrac{x}{\sqrt{\l}}=re^{-i\vt}\in S[-\tfrac{\pi}{2},0],
  \quad \l=|\l|e^{2i\vt}
\end{equation}
reduces the study of the family of curves
$\Gamma_\l=z_\l(\mathbb{R}_+)$,
$\l\in\overline{\mathbb{C}}_+$ to the study of the family
$e^{2i\vt}\xi(e^{-i\vt}\mathbb{R}_+)$ for $\vt\in[0,\tfrac{\pi}{2}]$.
We use the  representation of $t\in e^{-i\vt}\mathbb{R}_+$ in the
form
\begin{equation}
\label{AngularXi}
t=re^{-i\vartheta}=1+\eta e^{-i\varphi}, \qquad
\varphi\in[0,\pi], \quad r,\eta\ge0.
\end{equation}
 Using
this representation and (\ref{DefXi}), we have
\begin{equation}
\label{ChVarXi}
\xi(t)
=e^{-i\frac{3\varphi}{2}}\int_0^\eta
\sqrt{2+se^{-i\varphi}}\cdot s^{\frac{1}{2}}\,ds,
\quad
\partial_r\xi(t)
=e^{-i\vt}\sqrt{t^2-1}
=\sqrt{\eta}e^{-i(\vt+\frac{\varphi}{2})}\sqrt{1+re^{-i\vt}}.
\end{equation}

\begin{lemma} \label{Geom-0}
Fix $\vt\in[0,\frac{\pi}{2}]$ and let
$t\in e^{-i\vt}\mathbb{R}_+$. Write $t$ in the form (\ref{AngularXi}).
Then
\begin{enumerate}
\item \label{A7-3}
$\frac{2}{3}\sin^{\frac{3}{2}}\vt\le |\xi(t)|$;
if $|t-1|\le1$, then
$
\frac{2}{3} |t-1|^{\frac{3}{2}}\le |\xi(t)|\le  2
|t-1|^{\frac{3}{2}}
$,

\item \label{A7-2}
$
\arg \xi(t)\in
\left[
-\frac{3\varphi}{2}-\frac{\vartheta}{2},-\frac{3\varphi}{2}
\right]$,
$
-\varphi\in
[\frac{2}{3}\arg \xi(t),\frac{2}{3}\arg \xi(t)+\frac{\vartheta}{3} ]
$,

\item\label{A7.2-1}
$\arg{\partial_r}\xi(t)
\in [-\frac{\pi}{2}-\vt,-2\vt)
\cap
(-\frac{\varphi}{2}-\frac{3\vt}{2},-\frac{\varphi}{2}-\vt]$,

\item \label{A7-4}
if $\vt\in(0,\frac{\pi}{2}]$, then
$
\arg \xi(e^{-i\vt}\mathbb{R}_+)=[-\frac{3\pi}{2},-2\vt)$;
$
\arg \xi(e^{-i0}\mathbb{R}_+)=\{-\frac{3\pi}{2},0\}$,

\item \label{A7-5}
 if
$-\pi-\vt\le\arg \xi(t)
$,
then
$\arg\left( e^{2i\vt}{\partial_r}
\xi(t)\right)\in[-\frac{\pi}{3}+\frac{\vt}{6},\frac{\vt}{2}]$,

\item \label{A7-6}
 if $-\frac{\pi}{2}-2\vt\le\arg\xi(t)$, then
 $\arg\left( e^{2i\vt}{\partial_r}
\xi(t)\right)\in[-\frac{\pi}{6}-\frac{\vt}{6},\frac{\vt}{2}]$.
\end{enumerate}
\end{lemma}

{\it Proof.}
\ref{A7-3}. Assume $|t-1|\le1$ and consider the integrand in (\ref{ChVarXi}).
For $s\in[0,\eta]$ we have
$
1<\re\sqrt{2+se^{-i\varphi}}
$
and
$
|\sqrt{2+se^{-i\varphi}}|<3
$.
Substituting these estimates in (\ref{ChVarXi})
yields
$
\frac{2}{3}\eta^{\frac{3}{2}}\le |\xi(t)|
<3\cdot\frac{2}{3}\eta^{\frac{3}{2}}$,
$
\eta=|t-1|
$.
The relation
$\min\limits_{\arg t=-\vt}|t-1|=\sin\vartheta$ for
$\vartheta\in[0,\frac{\pi}{2}]$ finishes the proof.

 \ref{A7-2}. Consider the integrand in  (\ref{ChVarXi}).
 For $s\in[0,\eta]$ we have
$\arg\sqrt{2+se^{-i\varphi}}\in[-\frac{\vartheta}{2}, 0]
$, hence
$
\arg \xi(t)+\frac{3\varphi}{2}=
\arg \left(\int_0^\eta
\sqrt{2+se^{-i\varphi}}\cdot s^{\frac{1}{2}}\,ds\right)\in
[-\frac{\vartheta}{2},0]
$.

\ref{A7.2-1}.
In the second identity in (\ref{ChVarXi}) we have
 $(t^2-1)\in S[-\pi,-2\vt)$, so that
 $\arg \sqrt{|t|^2e^{-2i\vt}-1}\in[-\frac{\pi}{2},-\vt)$
 and
 $\arg{\partial_r}\xi(t)
\in [-\frac{\pi}{2}-\vt,-2\vt)$. Using
 $\arg\sqrt{1+re^{-i\vt}}\in(-\frac{\vt}{2},0]$
 in the same formula, we obtain
$\arg{\partial_r}\xi(t)
\in
(-\frac{\varphi}{2}-\frac{3\vt}{2},-\frac{\varphi}{2}-\vt]$.

\ref{A7-4}. Fix $\vt\in(0,\frac{\pi}{2}]$. By the
principle of boundary correspondence for
conformal mappings,
 $\xi(S[-\frac{\pi}{2},0])\subset S[-\frac{3\pi}{2},0]$,
 so $
\arg \xi(e^{-i\vt}\mathbb{R}_+)\subset [-\frac{3\pi}{2},0]$. Using
$
\xi(t)=t^2/2-\frac{1}{2}\ln t+O(1)
$,
 we have  $\arg \left(e^{2i\vartheta}\xi(t)\right)\to 0$ as
$r\to\infty$; by direct calculation, $\xi(0)=i\frac{\pi}{4}$.
Therefore,
$[-\frac{3\pi}{2},-2\vt)
\subset
\arg \xi(e^{-i\vt}\mathbb{R}_+)$. It remains to prove that
$\arg
\left(e^{2i\vartheta}\xi(t)\right)<0$.

By \ref{A7.2-1}, for a fixed $\vt$ the function
  $\Im \left(e^{2i\vartheta}\xi(re^{-i\vt})\right)$
    strictly decreases in $r$.
Therefore, using  $\xi( S[-\frac{\pi}{2},0])\subset S[-\frac{3\pi}{2},0]$
and $\xi(0)=i\frac{\pi}{4}$, we conclude that
the curve $ e^{2i\vt}\xi(e^{-i\vt}\mathbb{R}_+)$
crosses only the negative half of the imaginary axis.
Again using monotone decrease of
$\Im \left(e^{2i\vartheta}\xi(re^{-i\vt})\right)$ in $r$,
we obtain $\arg
\left(e^{2i\vartheta}\xi(t)\right)<0$,
which comletes the proof.

\ref{A7-5}.
Using the second identity in (\ref{ChVarXi}) and
$\sqrt{1+re^{-i\vt}}\in S[-\frac{\vartheta}{2},0]$,
we obtain
\begin{equation}
\label{Geom-5}
\arg
\left(
e^{2i\vartheta} {\partial_r}
\xi(t)
\right)
\in[-\tfrac{\varphi}{2}+\tfrac{\vartheta}{2},-\tfrac{\varphi}{2}+\vartheta].
\end{equation}
By hypothesis and \ref{A7-4},  $\arg \xi(t)\in[-\pi-\theta, -2\vartheta]$.
Thus using \ref{A7-2}  we obtain
$-\varphi\in[-\frac{2}{3}(\pi+\vartheta),
-\vartheta]$.
Substituting this into (\ref{Geom-5})
 proves \ref{A7-5}.

 \ref{A7-6}.
By \ref{A7-2},
$-\frac{\pi}{6}-\frac{2\vt}{3}
\le-\frac{\varphi}{2}\le0$.
Using $\arg \sqrt{re^{-i\vt}+1}\in[-\frac{\vt}{2},0]$, we obtain
$\arg (e^{2i\vartheta}{\partial_r}
\xi(t))=\arg(e^{i\vt}\sqrt{t^2-1})=
\arg(e^{i(\vt-\frac{\varphi}{2})}\sqrt{re^{-i\vt}+1})
\in[-\frac{\pi}{6}-\frac{\vt}{6},\frac{\vt}{2}]$.
$\blacksquare$

\begin{lemma}
\label{23}
Let $t=re^{-i\vt}$, $r\ge0$ for a fixed
$\vt\in[0,\frac{\pi}{2}]$. Then
\begin{enumerate}
\item \label{A7.2-2}
 if $\vt\in(0,\frac{\pi}{2}]$ and $r\ge0$
 then
$\Im\left( e^{2i\vt}{\partial_r}\xi(t) \right)<0$
and
$\Re\left( e^{2i\vt}{\partial_r}\xi(t) \right)>0$,

\noindent
if $\vt=0$ and  $r\in[0,1)$, then
$\Im\left( e^{2i\vt}{\partial_r}\xi(t) \right)<0$
and $\Re\xi(t)=0$,

\noindent
if $\vt=0$  and $r\in(1,\infty)$, then $\Im\xi(t)=0$ and
$\Re\left( e^{2i\vt}{\partial_r}\xi(t) \right)>0$,

\item  \label{A7.2-3}
if
$\arg\left( e^{2i\vt}
\xi(t)\right)\in (-\pi,-\frac{\pi}{2}+\vt)$, then
$
{\partial_r}\arg \xi(t)>0$,

\item \label{A.7.2-3a}
if $\arg\left( e^{2i\vt}
\xi(t)\right)\in[-\frac{3\pi}{2}+2\vt, -\frac{\vt}{4}]$,
then ${\partial_r}\arg \xi(t)\ge 0$,

\item  \label{A7.2-4}
if
$\arg\left( e^{2i\vt}
\xi(t)\right)\in (-\frac{3\pi}{2},-\pi+\vt)$, then
${\partial_r}|\xi(t)|<0$,\\
%
if
$\arg\left( e^{2i\vt}
\xi(t)\right)\in (-\frac{\pi}{2},\vt)$, then
${\partial_r}|\xi(t)|>0$.

\end{enumerate}
\end{lemma}
%

{\it Proof.}
\ref{A7.2-2}. For $\vt\in(0,\frac{\pi}{2}]$ the result follows from
 \refLem{A7.2-1}{Geom-0}.
 For $\vt=0$, the result follows from (\ref{DefXi}) by direct calculation.

\ref{A7.2-3}. Direct calculation yields
$
{\partial_r}\arg \xi(t)=
\frac{\partial}{\partial r}\Im \ln\xi(re^{-i\vt})=
\frac{|\xi'(t)|}{|\xi (t)|}
\sin\left\{ \arg{\partial_r}\xi (t)
-\arg \xi(t) \right\}
$.
By \refLem{A7.2-1}{Geom-0}, ${\partial_r}\arg \xi(t)$ is strictly positive for
$\arg \xi(t)\in(-\pi-2\vt,-\frac{\pi}{2}-\vt)$.

\ref{A.7.2-3a}. We have
$
\frac{\partial}{\partial r}\arg\xi(re^{-i\vt})
=
\Im\left(
\frac{{\partial_r}\xi(t)}{\xi(t)}
\right)
$. By \refLem{A7.2-1}{Geom-0},
$\arg {\partial_r}\xi(t)\in
[-\frac{\varphi}{2}-\frac{3\vt}{2},-\frac{\varphi}{2}-\vt]$, so
using \refLem{A7-2}{Geom-0} yields
$\arg ({\partial_r}\xi(t))\in
[\frac{1}{3}\arg\xi(t)-\frac{3\vt}{2},\frac{1}{3}\arg\xi(t)-\frac{5\vt}{6} ]
$. Hence,
$\arg \left(\frac{{\partial_r}\xi(t)}{\xi(t)}\right)
\in
[-\frac{2}{3}\arg\xi(t)-\frac{3\vt}{2},-\frac{2}{3}\arg\xi(t)-\frac{5\vt}{6 }]
$. By hypothesis and \refLem{A7-4}{Geom-0},
$\arg\xi(t)\in[-\frac{3\pi}{2},-2\vt-\frac{\vt}{4}]$,
which implies
$\arg \left(\frac{{\partial_r}\xi(t)}{\xi(t)}\right)
\in[0,\pi]$, as required.

\ref{A7.2-4}. We have
$
{\partial_r}|\xi(t)|^2=2 |\xi '(t)||\xi (t)| \cos\left\{
\arg{\partial_r}\xi (t)
-\arg \xi(t) \right\}
$.
By \refLem{A7.2-1}{Geom-0}, ${\partial_r}|\xi(t)|$ is positive for
$\arg \xi(t)\in (-\frac{\pi}{2}-2\vt,-\vt)$
 and negative for
 $\arg \xi(t)\in (-\frac{3\pi}{2}-2\vt,-\pi-\vt)$.
$\blacksquare$

\begin{lemma}
\label{DivG2}
 For each
 $\l=|\l|e^{2i\vt}\in
 \C_+\setminus\{0\}$
   there exists a unique $x_*\ge0$ such that
 $|z_\l(x_*)|=\min\limits_{x\ge0}|z_\l(x)|$.
 Define $t_*\in e^{-i\vt}\mathbb{R}_+$ and $z_*\in\Gamma_\l$
 by $t_*=\frac{x_*}{\sqrt{\l}}$ and
$z_*=z_\l(x_*)$. Then
\begin{enumerate}
\item \label{A7.3-2}
   $|z_\l(\cdot)|$ is strictly decreasing on $[0,x_*)$ and
 strictly  increasing on $(x_*,\infty)$,
\item \label{27}
for $\l$ fixed $\Re z_\l(\cdot)^{\frac{3}{2}}$ is  non-decreasing
on $\mathbb{R}_+$. If $\l>0$ and $x\in[x_*,\infty)$, or if $0<\arg\l\le\pi$
 and
 $x\in\mathbb{R}_+$, then it is strictly increasing,

\item \label{A7.3-1}
if $\vt=0$, then $z_*=0$ and $x_*=\sqrt{\l}$,\\
if $\vt\in(0,\frac{\pi}{2}]$, then
$\arg z_*\in
[-\frac{\pi}{2}-\frac{\vartheta}{2},-\frac{\pi}{2}+\frac{5}{6}\vt]$,
$-\pi+\frac{\pi}{22}\le\arg \left( e^{2i\vt}\xi(t_*)\right)$,
    \item \label{A7.3-3}
   if $t=e^{-i\vt}r$
  for  $r\in [0,|t_*|)$, then
  $\arg\left( e^{2i\vt}{\partial_r}\xi(t)\right)
  \in[-\frac{\pi}{2}+\vt,-\frac{\pi}{4}+\frac{3\vt}{4}]$,
\item \label{A7.3-4}
$|t_*|\le\sqrt{2}$.

\end{enumerate}
\end{lemma}

{\it Proof.}
First we prove existence and uniqueness of $x_*$.
 By (\ref{Z3-2PowInXi}), it is sufficient
to prove that there exist a unique minimum of $|\xi(t)|$
on $e^{-i\vt}\mathbb{R}_+$.
For $\arg\l=0$ the result is evident.
Fix $\vt\in(0,\frac{\pi}{2}]$ and let $t=re^{-i\vt}$, $r\ge0$.
Direct calculation yields $\partial_r|\xi(t)|^2|_{r=0}=-(\pi/\sqrt{2})\cos\vt\le0$.
By \refLem{A7-4}{Geom-0} and \refLem{A7.2-4}{23}, we have
$\partial_r|\xi(t)|>0$ as $r\to\infty$.  Therefore
 there exists at least one
point $t\in e^{-i\vt}\mathbb{R}_+$ such that ${\partial_r}|\xi(t)|=0$.
This point is unique if
\begin{equation}
\label{DerXi}
\tfrac{1}{2}\tfrac{\partial^2}{\partial r^2}|\xi(t)|^2
=\left|\xi'(t)\right|^2+
\Re\left(e^{-2i\vt}\xi''(t)\overline{\xi(t)}\right)>0
\qquad\text{for}\quad
 t\in S[-\tfrac{\pi}{2},0]
.
\end{equation}
By
$
\xi'(t)=\sqrt{t^2-1}$
and
$
\xi''(t)=\frac{t}{\sqrt{t^2-1}}
$,
it is sufficient to show that
\begin{equation}\label{XiW}
 \left|\frac{\xi(t)}{w(t)}\right|<1
 \qquad \hbox{for}\quad \Re t\ge0,
\end{equation}
where
$w(t)=\frac{(\xi'(t))^2}{\xi''(t)}=\frac{(t^2-1)^{\frac{3}{2}}}{t}$ on
$\mathbb{C}\setminus(-\infty,1]$ and $w(t)>0$ for $t>1$.
Note that  $\frac{\xi(t)}{w(t)}$ is analytic in the
half-plane $\Re t>-1$, since  both $w$ and $\xi$
change sign when crossing the cut
 $(-1,1]$. Thus we prove (\ref{XiW})
applying the  principle of maximum for the expanding
half-disks $D_R=\{t: \Re t\ge0, |t|\le R\}$. For $t\to\infty$
we have $w(t)\sim t^2$
and $\xi(t)\sim \frac{t^2}{2}$  uniformly in $|\arg t|\le\frac{ 3\pi}{4}$.
Hence $|\frac{\xi(t)}{w(t)}|\to\frac{1}{2}$
uniformly on the arcs
$|t|=R$, $|\arg t|\le\frac{\pi}{2}$ as $R\to\infty$. For $\Re t=0$ deformation of
the integration path in (\ref{DefXi}) gives
$
|\xi(t)|^2\le (\frac{\pi}{4})^2+|t|^2(1+\frac{|t|}{2})^2
$,
so that $|\xi(t)|\le\frac{7}{12}|w(t)|$.
Thus the inequality
$|\frac{\xi}{w}|<\frac{7}{12}$ holds
on the boundary of
$D_R$ for sufficiently large $R$.
By the maximum principle, this yields (\ref{XiW}), (\ref{DerXi})
and uniqueness of the point $t_*\in e^{-i\vt}\mathbb{R}_+$
satisfying ${\partial_r}|\xi(t_*)|=0$.
By (\ref{Z3-2PowInXi}), $x_*=\frac{t_*}{\sqrt{\l}}$ is the unique
minimum of $|z_\l(\cdot)|$ on $\mathbb{R}_+$, as required.

\ref{A7.3-2}. The result follows from (\ref{Z3-2PowInXi}), (\ref{DerXi}) and
$\partial_r|\xi(t)|^2|_{r=0}=-(\pi/\sqrt{2})\cos\vt\le0$.

\ref{27}. By (\ref{Z3-2PowInXi}), for
$\frac{x}{\sqrt{\l}}=re^{-i\vt}$
the sign of $\partial_x \Re z_\l(x)^{\frac{3}{2}}$
is the same as that of $\partial_r\Re(e^{2i\vt}\xi(re^{-i\vt}))$.
Therefore, the result follows from \refLem{A7.2-2}{23}.

\ref{A7.3-1}. For $\arg\l=0$ the result is evident; fix
$\vt\in(0,\frac{\pi}{2}]$ and let $t=re^{-i\vt}$, $r\ge0$.
We have
${\partial_r}|\xi(t)|^2=
2\Re
\left(
\xi(t)\overline{\partial_r \xi(t)}
\right)
$.
By \refLem{A7-2}{Geom-0}, for
$t=1+\eta e^{-i\varphi}\in S[-\frac{\pi}{2},0) $  we have
$\arg\xi(t)\in
[-\frac{3\varphi}{2}-\frac{\vartheta}{2}, -\frac{3\varphi}{2}]$;
by \refLem{A7.2-1}{Geom-0}, we have
$\arg\overline{\partial_r\xi(t)} \in
[\frac{\varphi}{2}+\vt,\frac{\varphi}{2}+\frac{3\vartheta}{2}]$.
Thus
$\arg\left(
\xi(t)\overline{\partial_r \xi(t)}
\right) \in
[-\varphi+\frac{\vartheta}{2}, -\varphi+\frac{3\vt}{2}]$
and we have
$
{\partial_r}|\xi(t)|< 0$ for $-\varphi\in
\left(-\frac{3\pi}{2}-\frac{\vartheta}{2},
-\frac{\pi}{2}-\frac{3\vt}{2}\right)$
and
$
{\partial_r}|\xi(t)|> 0 $ for $-\varphi\in
\left(-\frac{\pi}{2}-\frac{\vartheta}{2},
\frac{\pi}{2}-\frac{3\vt}{2}\right)$.
Therefore,
the point $t_*=r_*e^{-i\vt}=1+\eta_* e^{-i\varphi_*}$
satisfies
\begin{equation}
\label{Ab.eq}
-{\varphi_*}\in \text{$\left[-\tfrac{\pi}{2}-\tfrac{3{\vartheta}}{2},
-\tfrac{\pi}{2}-\tfrac{{\vartheta}}{2}\right]$}.
\end{equation}
 Using (\ref{Ab.eq}) and
\refLem{A7-2}{Geom-0} we obtain
  \begin{equation}\label{Argxi(z*)}
\arg\xi(t_*)\in
\left[-\tfrac{3\pi}{4}-\tfrac{11\vartheta}{4},
-\tfrac{3\pi}{4}-\tfrac{3\vartheta}{4}
\right].
\end{equation}
By (\ref{Z3-2PowInXi}), this  proves
$\arg z_*\in
[-\frac{\pi}{2}-\frac{\vartheta}{2},-\frac{\pi}{2}+\frac{5}{6}\vt]$.
By (\ref{Argxi(z*)}) for $0\le\vt\le\frac{3\pi}{11}$
and  \refLem{A7-4}{Geom-0} for $\frac{3\pi}{11}\le\vt\le\frac{\pi}{2}$,
we have $-\pi+\frac{\pi}{22}\le\arg(e^{2i\vt}\xi(t_*))$, as required.

\ref{A7.3-3}. By (\ref{Ab.eq}), for
$t=re^{-i\vt}=1+\eta e^{-i\varphi}$ and
 $r<|t_*|$
 we have
$
   -\pi\le-\varphi\le-\frac{\pi}{2}-\frac{\vartheta}{2}
$.
Using \refLem{A7.2-1}{Geom-0} we obtain
$
-\frac{\pi}{2}+\vt\le
\arg\left(
e^{2i\vt}{\partial_r}\xi(t)
\right)
\le \vt-\frac{\varphi}{2}\le-\frac{\pi}{4}+\frac{3\vt}{4}.
$

\ref{A7.3-4}. For $\vt=0$ the result is evident.
For a fixed $\vt\in(0,\frac{\pi}{2}]$ write
$t=re^{-i\vt}=1+\eta e^{-i\varphi}$.
Geometric considerations show that
$
r=\frac{\sin\varphi}{\sin(\varphi-\vt)}
$.
By (\ref{Ab.eq}),
 for $t_*=r_*e^{-i\vt}=1+\eta_* e^{-i\varphi_*}$ we have
$
-\frac{\pi}{2}-\frac{\vt}{2}\le\vt-\varphi_*\le-\frac{\pi}{2}+\frac{\vt}{2},
$
so that
$r_*=\frac{\sin\varphi_*}{\sin(\varphi_*-\vt)}\le\frac{1}{\sin\frac{\pi}{4}}=
\sqrt{2}$.
$\blacksquare$

In the next Lemma we analyse the  curves
 $\Gamma_\l^{\pm}$, defined in (\ref{DefGPM}).

\begin{lemma}
\label{SvaGForAll}
Let
$\l=|\l|e^{2i\vt}\in\ol\C_+\setminus\{0\}$. Then
\begin{enumerate}
 \item \label{A7.5-1}
 if $\vt=0$, then
 $\Gamma_\l^-=[-(\frac{3\pi}{8}\l)^{\frac{2}{3}},0]$,
$\Gamma_\l^+=[0,+\infty)$,
 \\
 if $\vt\in(0,\frac{\pi}{2}]$, then $\Gamma_\l\subset
 S[-\pi+\frac{4\vt}{3},0)$,
 $\Gamma_\l^-\subset S[-\pi+{\frac{4\vt}{3}},-\frac{\pi}{2}+\frac{5\vt}{6}]$,
$\Gamma_\l^+\subset S[-\frac{\pi}{2}-\frac{\vartheta}{2},0)$,
   \item \label{A7.5-4} if $0<\arg\l\le\delta$, then
   $\Gamma_\l^-\subset S[-\pi+0,-\frac{\pi}{3}-\frac{\pi}{12}]$,
   $\Gamma_\l^+\subset S[-\frac{\pi}{2}-\frac{\pi}{20},0)$,
 \item \label{A7.5-2}
 $\inf\limits_{z\in\Gamma_\l}|z|\ge |\l|^{\frac{2}{3}}\sin\vt$,
 \item \label{A7.5-3}
 $\Gamma_\l^-\subset\{z: |z|
 \le
\left(\frac{3\pi}{8}|\l|\right)^{\frac{2}{3}}
\}$,
  the length of $\Gamma_\l^-$ satisfies
  $|\Gamma_\l^-|\le \co|\l|^{\frac{2}{3}}$.
\end{enumerate}

\end{lemma}
{\it Proof  } \ref{A7.5-1}. For $\vt=0$ the result is evident.
 For $\vt\in(0,\frac{\pi}{2}]$  the assertion on $\Gamma_\l$ follows from
 \refLem{A7-4}{Geom-0} and relation (\ref{Z3-2PowInXi}).

Now prove  the assertion on $\Gamma_\l^\pm$.
 By \refLem{A7.3-1}{DivG2}, we have only show that
$\arg z_\l(x)$ is non-decreasing at $x=x_*$.
Writing \refLem{A7.3-1}{DivG2} in terms of
$t_*=r_*e^{-i\vt}=1+\eta_* e^{-i\varphi_*}$ yields
 (\ref{Argxi(z*)}).  Hence the hypothesis of \refLem{A.7.2-3a}{23}
 is satisfied for $t_*$, so that ${\partial_r}\arg \xi(t_*)\ge0$.
 By (\ref{Z3-2PowInXi}),
 $\partial_x\arg z_\l(x_*)\ge0$. This
 completes the proof.

\ref{A7.5-4}. By \ref{A7.5-1} of this Lemma, we have $\Gamma_\l^-
\subset S[-\pi, -\frac{\pi}{2}+\frac{5\delta}{12}]$
and
$\Gamma_\l^+
\subset S[-\frac{\pi}{2}-\frac{\delta}{4},0)$. By definition (\ref{DefDelta}),
$\delta<\frac{\pi}{5}$,
which implies the result.

\ref{A7.5-2}. The result follows from \refLem{A7-3}{Geom-0} and
relation (\ref{Z3-2PowInXi}).

\ref{A7.5-3}.
The inclusion follows from \refLem{A7.3-2}{DivG2} and (\ref{Defz0zE}).
Using the definition  (\ref{DefZofX}), we obtain
$|\Gamma_\l^-|=|\l|^{\frac{2}{3}}\int_0^{|t_*|}|\eta'(re^{-i\vt})| dr$.
The integral is  bounded uniformly in
 $\l$, since
$\eta(t)$ is analytic for $\Re t>-1$
and $|t_*|\le\sqrt{2}$  (see
 \refLem{A7.3-4}{DivG2}). This proves
the estimate for $|\Gamma_\l^-|$.
$\blacksquare$

\section{The main estimates}
\label{Section:FinalEstimates} In this section we prove the main
estimates (\ref{EstPsi-0}--\ref{EstPsi-*}). We  consider only the
case $\Im\l\ge0$,
 since for $\Im\l\le0$ the results are analogous.
In order to formulate the result
we recall the necessary definitions and notations.

By
 (\ref{DefZofX}), we have
$ z_\l(x)= \left(
\frac{3}{2}\l\xi(\frac{x}{\sqrt{\l}})\right)^{\frac{2}{3}} $, $
\xi(t)= \int_1^t\sqrt{s^2-1}\, ds $, where $\xi(t)>0$ for $t>1$. By
Definition~\ref{Z-star}, $|x_\l(\cdot)|$ on $\mathbb{R}_+$ has the
unique minimum in $x_*\equiv x_*(\l)$; for $\l>0$ $x_*=\sqrt{\l}$ is
the turning point.
 By (\ref{Example}),
$
\phi(\l)=
2^{\frac{3}{4}}\sqrt{\pi}
\left(\frac{\l}{2e}\right)^{\frac{\l}{4}}
$ takes its principal value on $\mathbb{C}\setminus\mathbb{R}_-$.
According to (\ref{Example}),
 $$
\rho(x,\l)=
\frac{
1+|\l|^{\frac{1}{2}}+|x^2-\l|^{\frac{1}{2}}
}{
1+|\l|^{\frac{5}{12}}+|x^2-\l|^{\frac{5}{4}}
}.
 $$

For $\arg\l>0$ we define by $\sqrt{x^2-\l}$ the branch,
 analytic on $\mathbb{R}_+$,
 such that $\arg\sqrt{x^2-\l}\to0$ as $x\to+\infty$.
For $\arg\l=0$ we set $\sqrt{x^2-\l}=\sqrt{x^2-(\l+i0)}$.

Using Definition~\ref{DefHDelta} for $H_{\pm\delta}(\l)$ and (\ref{DefD_Z^Epsilon})
for $D_Z^\delta(\l)$, introduce the domains
\begin{align*}
D_0^\delta(\l)
    &\overset{\text{\rm def}}{=}
D_Z^\delta(\l)\cap\{ z:|\arg z|\le\pi-\tfrac{\delta}{3}\},
    &\l  &\in S_{1/2}[0,\pi],\\
D_{+}^\delta(\l)
    &\overset{\text{\rm def}}{=}
     D_Z^\delta(\l)\cap
     \left(
     S[\tfrac{\pi}{3}+\tfrac{4\vt}{3}+\tfrac{\delta}{3},
\pi-\tfrac{\delta}{3}]
\cup
H_{-\delta}(\l)
\cup
S[-\pi+\tfrac{4\vt}{3}, \tfrac{\pi}{3}-\tfrac{\delta}{3}]
\right),
    & \l &\in S_{1/2}[0,\pi-\delta],\\
D_{-}^\delta(\l)
    &\overset{\text{\rm def}}{=}
 D_Z^\delta(\l)
 \cap
 \{z:|\tfrac{\pi}{3}-\arg z|\ge\tfrac{\delta}{3}\},
    &\l &\in S_{1/2}[0,\pi],\\
D_*^\delta(\l)
    &\overset{\text{\rm def}}{=}
    D_Z^\delta(\l)\cap
\left(
    S[\tfrac{\pi}{3}+\tfrac{\delta}{3},
\pi+\tfrac{4\vt}{3}]
\cup
H_{\delta}(\l)
\cup
S[-\tfrac{\pi}{3}+\tfrac{\delta}{3},
 -\tfrac{\pi}{3}+\tfrac{\delta}{3}+\tfrac{4\vt}{3}]
\right),
    & \l &\in S_{1/2}[\delta,\pi],
\end{align*}
\begin{figure}[h!]
\centering
\includegraphics[width=165mm]{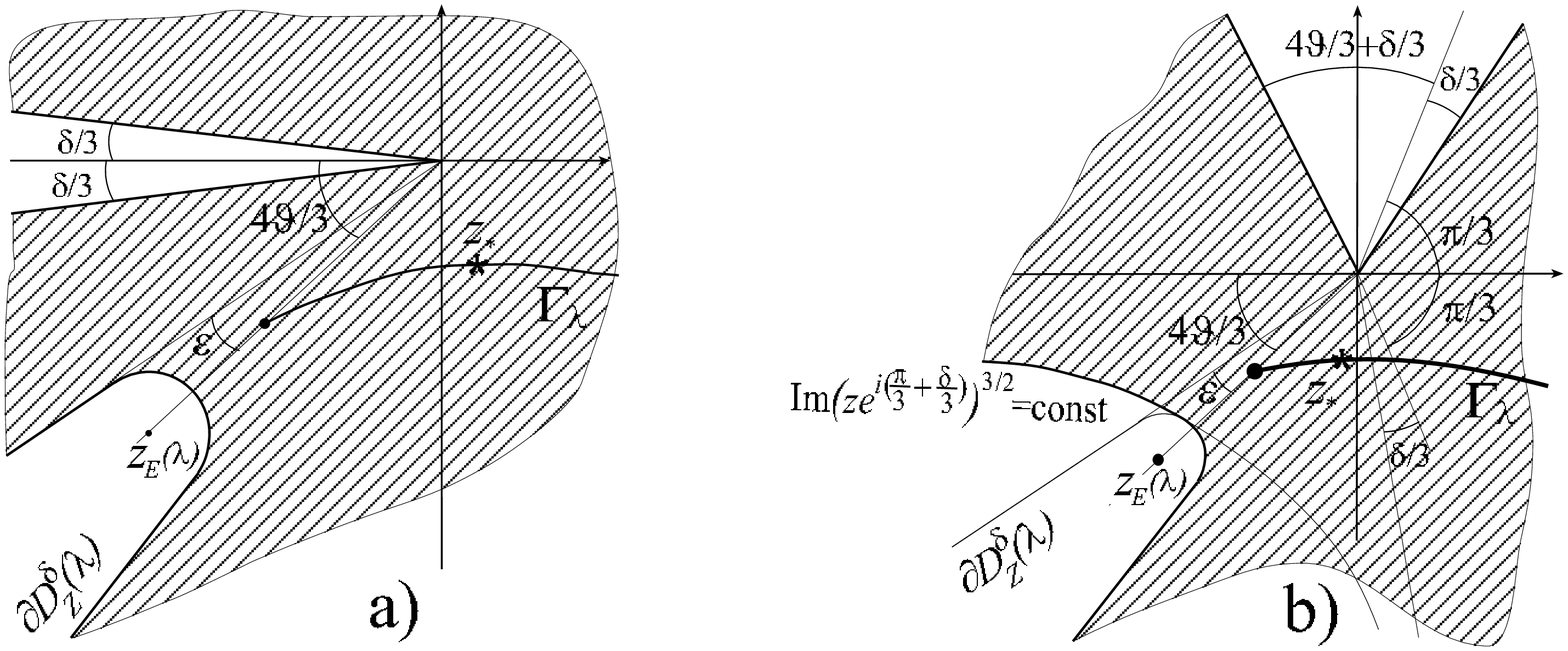}\\
\vskip2em
\includegraphics[width=165mm]{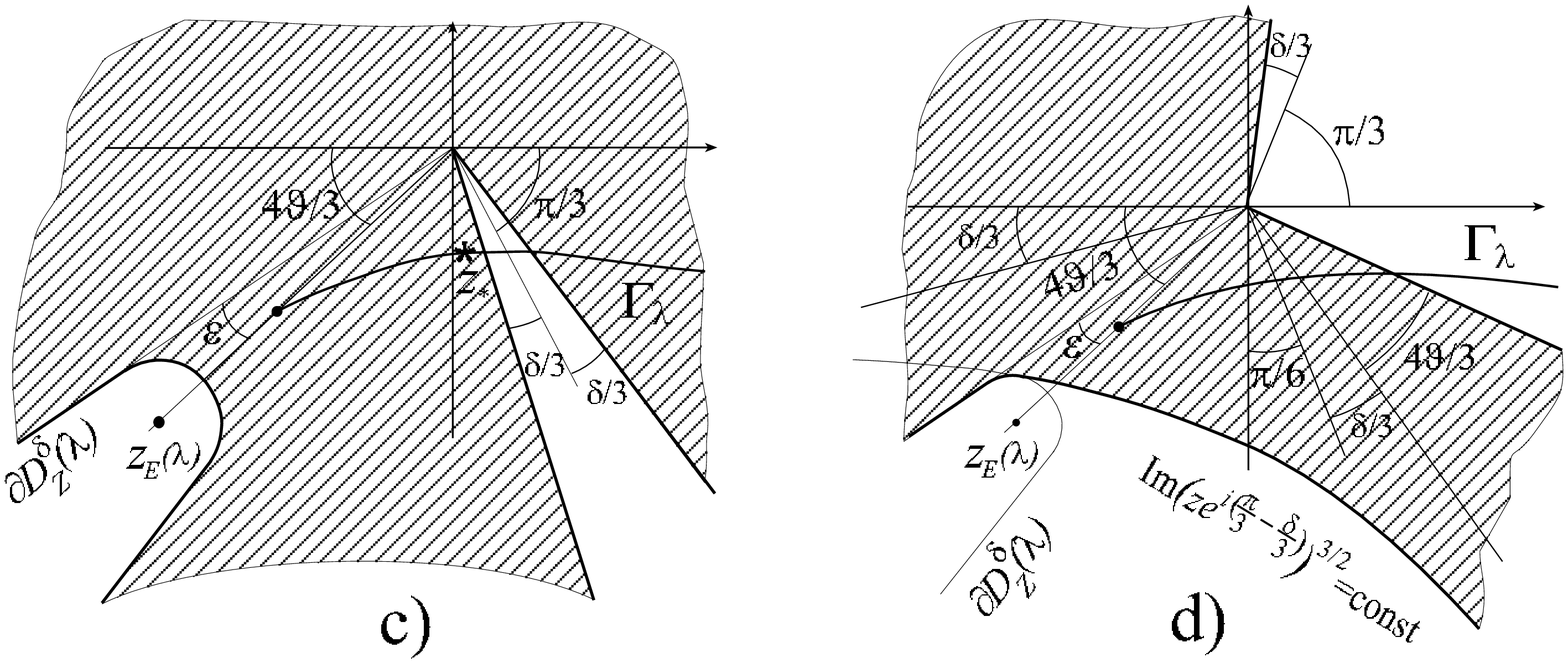}
\caption[short caption here]{ The domains a) $D_0^\delta(\l)$,
b) $D_+^\delta(\l)$,
c) $D_-^\delta(\l)$,
d) $D_*^\delta(\l)$,
 where uniform estimates
(\ref{MainEstimateForASmall})
are fulfilled.
Here $2\vt=\arg\l$.}
\label{EstimatesOfAandA+}  
\end{figure}
which are schematically presented of Fig.\ref{EstimatesOfAandA+}.

 We are interested in  the range
 of $\arg\l$ such that $\Gamma_\l^\pm$, given by (\ref{DefGPM}),
  are within the estimate's validity domain.
By \refLem{A7.2-1}{Geom-0},
$\Im (z_\l(x)e^{i(\frac{\pi}{3}-\frac{\delta}{3})})^{\frac{3}{2}}$ is non-decreasing
in $x$ for $\delta\le\arg\l\le\pi$. Taking into account the definition
(\ref{DefD_Z^Epsilon}), we conclude that
 for $\delta\le\arg\l\le\pi$
  the curve
$\Im (ze^{i(\frac{\pi}{3}-\frac{\delta}{3})})^{\frac{3}{2}}=${\rm const},
tangential to the boundary of $D_Z^\delta(\l)$, does not intersect $\Gamma_\l$,
\begin{equation}\label{Gamma_LambdaDoesNotIntersect}
\Gamma_\l\cap
\{
z\in S[-\pi+\tfrac{\delta}{3},-\tfrac{\pi}{3}+\tfrac{\delta}{3} ]:
\Im (ze^{i(\frac{\pi}{3}-\frac{\delta}{3})})^{\frac{3}{2}}
=
\Im (w_{\delta}(\l)
e^{i(\frac{\pi}{3}-\frac{\delta}{3})})^{\frac{3}{2}}
\}=\emptyset,
\qquad
\l\in S_{1/2}[\delta,\pi].
\end{equation}

\begin{theorem}
\label{Basic:InVariableX}
 Let $0<\delta<\frac{\pi}{5}$ and let $\psi$ be the solution of
 Eq.(\ref{OurEq}), satisfying (\ref{DefinitionParabolic CylFunct}).
  Then there exist a positive number
$C_\delta$, independent of $x$ and $\l$, such that

 For $0\le\arg \l\le\pi$ and $z_\l(x)\in D_0^\delta(\l)$ we have
\begin{equation}\label{EstPsi-0}
\left|
\psi'(x,\l)
+\psi(x,\l)\sqrt{x^2-\l}
\right|
\le C_\delta
|\phi(\l)
\rho(x,\l)
e^{-\l\xi(\tfrac{x}{\sqrt{\l}})}|.
\end{equation}
In particular, the estimate holds for
$\delta\le\arg \l\le\pi$, $x\in\mathbb{R}_+$ and
$0\le\arg \l\le\delta$, $x\in[x_*(\l),\infty]$.

For $0\le\arg \l\le\pi-\delta$ and
$z_\l(x)\in D_+^\delta(\l)$ we have
\begin{equation}\label{EstPsi-+}
\left|
\psi'(ix,-\l)
+i\psi(ix,-\l)\sqrt{x^2-\l}
\right|
\le C_\delta
|\phi(-\l)e^{-i\frac{\pi}{2}\l}
\rho(x,\l)
e^{\l\xi(\tfrac{x}{\sqrt{\l}})}|.
\end{equation}
In particular, the estimate holds for
$0\le\arg \l\le\pi-\delta$ and $x\in\mathbb{R}_+$.

For $0\le\arg \l\le\pi$ and
$z_\l(x)\in D_-^\delta(\l)$ we have
\begin{equation}\label{EstPsi--}
\left|
\psi'(-ix,-\l)
+i\psi(-ix,-\l)\sqrt{x^2-\l}
\right|
\le C_\delta
|\phi(-\l)
\rho(x,\l)
e^{-\l\xi(\tfrac{x}{\sqrt{\l}})}|.
\end{equation}
In particular, the estimate holds for
 $0\le\arg \l\le\delta$ and $x\in[0, x_*(\l)]$.

For $\delta\le\arg \l\le\pi$ and
$z_\l(x)\in D_*^\delta(\l)$ we have
\begin{equation}\label{EstPsi-*}
\left|
\psi'(-x,\l)
+\psi(-x,\l)\sqrt{x^2-\l}
\right|
\le C_\delta
|\phi(\l)e^{-i\frac{\pi}{2}\l}
\rho(x,\l)
e^{\l\xi(\tfrac{x}{\sqrt{\l}})}|.
\end{equation}
In particular, the estimate holds for
$\frac{\pi}{2}-\frac{\delta}{2}\le\arg \l\le\pi$ and $x\in\mathbb{R}_+$.
\end{theorem}

The plan of the proof is as follows.
In Definition~\ref{DefinitionOfA}  we introduce the
solutions $A_0$, $A_\pm$ and $A_*$ of the parabolic cylinder equation
in $z$-coordinate (\ref{EqinZ}).
Each solution is asymptotically close to  one of
the Airy functions $\Ai(z)$, $\Ai(ze^{\pm i\frac{2\pi}{3}})$ in some sector.
It is convenient to separate explicitly
the exponential multipliers of the solutions and
 proceed in terms of the modified
ones, which we denote by $ a_0$, $ a_\pm$ and $ a_*$
(see Definition~\ref{DefASmall}).
In Theorem~\ref{Basic:Svoistva}
 we  estimate the $z$-derivative of the modified solutions.
 For each modified solution we specify the domain where
  the estimate is uniform in both
 $z$ and $\l$; we find the range
 of $\arg\l$ such that $\Gamma_\l^+$ or $\Gamma_\l^-$ (or both)
 is within the domain.
  Finally we transform the estimates of
$\partial_z a_0$, $\partial_z a_\pm$ and $\partial_z a_*$
into the estimates of
 Theorem~\ref{Basic:InVariableX}
 using the relation of $A_\nu(z,\l)$
 to the parabolic cylinder
 functions $\psi(\pm x,\l)$, $\psi(\pm ix,-\l)$.
 We also present connection formulas and calculate Wronskians
  for solutions $A_0$, $A_\pm$ and $A_*$.


\begin{definition}
\label{DefinitionOfA}
$A$, $A_\pm$ and $A_*$ are the solutions of Eq.(\ref{EqinZ}),
satisfying the following asymptotics as $|z|\to\infty$:

\begin{equation}\label{DefA_0}
A_0(z,\l)=\Ai(z)(1+O(z^{-\frac{3}{2}})), \quad
\partial_zA_0(z,\l)=\Ai'(z)(1+O(z^{-\frac{3}{2}}))
\end{equation}
 for $\l\in S_{1/2}[0,\pi]$ and $z\in
S[-\pi+\frac{4\vt}{3}+\ve, \pi-\frac{\delta}{3}]$;
\begin{equation}\label{DefA_+}
A_+(z,\l)=\Ai(z\omega)(1+O(z^{-\frac{3}{2}})), \quad
\partial_zA_+(z,\l)=\omega\Ai'(z\omega)(1+O(z^{-\frac{3}{2}}))
\end{equation}
for $\l\in S_{1/2}[0,\pi-\delta]$ and
$z\in S[-\pi+\frac{4\vt}{3}+\ve, \frac{\pi}{3}-\frac{\delta}{3}]$;
\begin{equation}\label{DefA_-}
A_-(z,\l)=\Ai(z\overline{\omega})(1+O(z^{-\frac{3}{2}})), \quad
\partial_zA_-(z,\l)=
\overline{\omega}\Ai'(z\overline{\omega})(1+O(z^{-\frac{3}{2}}))
\end{equation}
for $\l\in S_{1/2}[0,\pi]$ and
$z\in S[-\frac{\pi}{3}+\frac{\delta}{3}, \pi+\frac{4\vt}{3}-\ve]$;
\begin{equation}\label{DefA_*}
A_*(z,\l)=\Ai(z\omega)(1+O(z^{-\frac{3}{2}})), \quad
\partial_zA_*(z,\l)=\omega\Ai'(z\omega)(1+O(z^{-\frac{3}{2}}))
\end{equation}
for $\l\in S_{1/2}[\delta,\pi]$ and
$z\in S[\frac{\pi}{3}+\frac{\delta}{3}, \pi+\frac{4\vt}{3}-\ve]$.
Here  $\delta$ and $\ve$ be given by (\ref{DefDelta}) and
(\ref{DefD_Z^Epsilon}), respectively.
\end{definition}

\begin{figure}[h!]
\centering
\includegraphics[height=6cm]{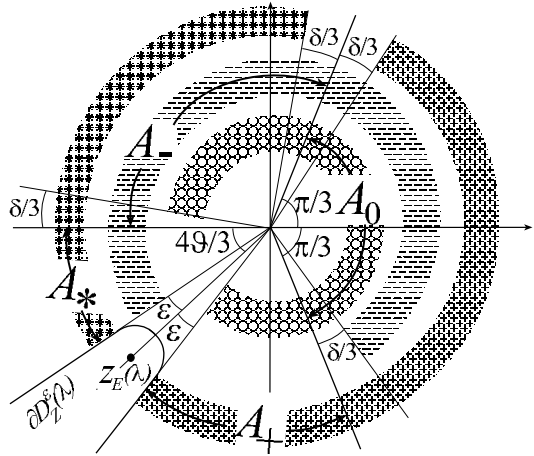}\\
\caption[short caption here]{ The sectors on $z$-plane,
 where $A_0$, $A_\pm$ and $A_*$
are asymptotically close  to Airy functions (here $2\vt=\arg\l\in[0,\pi]$).
Arrows indicate the sectors where the solutions decay. }
\label{AsymptSectorsAll}  
\end{figure}

The sectors from the definition are schematically
 presented on Fig.\ref{AsymptSectorsAll}.

The functions $A_0$ and $A_-$ are asymptotic to
$\Ai(z)$ and $\Ai(z\overline{\omega})$, respectively.
The sector of exponential decay
of $\Ai(z\omega)$  is divided in two subsectors
by the ray $\arg z=-\pi+\frac{4\vt}{3}$ (see (\ref{DefDomainDZ})).
The solution
 $A_+$ is asymptotic to $\Ai(z\omega)$
 in one  subsector, $A_*$ --- in another one.

Further analysis  of the perturbed Airy equation
(\ref{EqinZ}) use the following connection
and asymptotic formulas for the
Airy functions (see \cite{AS}):
 \begin{equation}
\label{BfromA} \{\Ai(z),\Bi(z)\}=1, \quad
 \Bi(z)=i\left(2e^{-i\frac{\pi}{3}}\Ai(\omega z)-\Ai(z)\right),\ \
 \omega=e^{i\frac{2\pi }{3}},
 \quad
\end{equation}
\begin{equation}
\label{ab1} \Ai(z)=e^{-i\frac{\pi}{3}}\Ai(z\omega)+
e^{i\frac{\pi}{3}}\Ai(z\ol\omega),\quad
\Bi(z)=ie^{-i\frac{\pi}{3}}\Ai(z\omega)-
ie^{i\frac{\pi}{3}}\Ai(z\ol\omega).
\end{equation}
\begin{equation}
\label{AsAiry}
\Ai(z)=\frac{e^{-\frac{2}{3}z^{\frac{3}{2}}}}{2\sqrt{\pi}\sqrt[4]{z}}
  \left(1+O(z^{-\frac{3}{2}})\right),\qquad
 |z|\to\infty,\qquad
 |\arg z|<\pi-\epsilon,
 \quad \forall\epsilon>0,
\end{equation}
Next we introduce
the modified solutions $a_\nu$, where the exponential
multipliers are separated explicitly. The basic
estimates are formulated in terms of these modified solutions.
\begin{definition}\label{DefASmall}
 For $|\l|>1/2$ and $z\in D_Z^0(\l)$
\begin{align}
a_0(z,\l)&=e^{\frac{2}{3}z^{\frac{3}{2}}} A(z,\l),
&
\quad\text{for}\quad \l &\in S(-\pi,\pi),
&
 z &\not\in \mathbb{R}_-,\\
\label{}
a_+(z,\l)&=e^{\frac{2}{3}(z\omega)^{\frac{3}{2}}} A_+(z,\l),
&
\quad\text{for}\quad \l &\in S(-\pi,\pi),
&
z &\not\in  e^{i\frac{\pi}{3}}\mathbb{R}_+,
\\
a_-(z,\l)&=e^{\frac{2}{3}(z\overline{\omega})^{\frac{3}{2}}}
A_-(z,\l),
&
\quad\text{for}\quad \l &\in S(-\pi,\pi),
&
 z &\not\in
e^{-i\frac{\pi}{3}}\mathbb{R}_+,
\\
\label{}
a_*(z,\l)&=e^{\frac{2}{3}(z\omega)^{\frac{3}{2}}} A_*(z,\l),
&
\quad\text{for}\quad \l &\in\mathbb{C}\setminus \mathbb{R},
&
\quad z &\not\in
 e^{i\frac{\pi}{3}}\mathbb{R}_+.
\end{align}
\end{definition}

The next theorem is essentially Theorem~\ref{Basic:InVariableX}
in $z$-variable. It
 gives the estimates of the $z$-derivatives of $a_\nu$,
 uniform in both $z$ and $\l$.
It is proved using  an equivalent integral equation.
The domains of the estimates
are defined in the beginning of this section.
 Their complicated structure is the cost of transition from asymptotic formulas
to uniform estimates.
\footnote{Derivation of the classical estimates
(\ref{OlversEstimates}) in \cite{Olver1960} also involves detailed discussion of domains.}

\begin{theorem}
\label{Basic:Svoistva}
  The solutions $A_0$, $A_\pm$ and $A_*$
exist and are defined uniquely. The solutions are
analytic in both $\l$ and $z$ for
$\l\in S_{1/2}(-\pi,\pi)$ and $z\in D_Z(\l)$. The asymptotics
(\ref{DefA_0}-- \ref{DefA_*}) are uniform in both $z$ and $\l$,
and the error term can be replaced by
$O(z^{-\frac{1}{2}}\l^{-\frac{2}{3}})$.

Fix $\delta\in(0,\frac{\pi}{5})$ and
let $D_Z^\delta(\l)$  be given by (\ref{DefD_Z^Epsilon}).
Then
the estimate
\begin{equation}\label{MainEstimateForASmall}
|\partial_z a_\nu(z,\l)|\le\frac{\co_\delta}{(1+|z|)^{\frac{5}{4}}},
\end{equation}
where $C_\delta$ is independent of $z$ and $\l$,
 is fulfilled  in the following cases:
\begin{description}
    \item $\nu=0$ and
$z\in D_0^\delta(\l)$,
$\l\in S_{1/2}[0,\pi]$.
In particular,
the  estimate  is valid
 for $\l\in S_{1/2}[\delta,\pi]$, $z\in\Gamma_\l$
 and for $\l\in S_{1/2}[0,\delta]$,  $z\in\Gamma_\l^+$.
    \item $\nu=+$ and
    $z\in D_+^\delta(\l)$,
 $\l\in S_{1/2}[0,\pi-\delta]$.
In particular,
the estimate is valid
 for $\l\in S_{1/2}[0,\pi-\delta]$, $z\in\Gamma_\l$.
    \item $\nu=-$ and
$z\in D_-^\delta(\l)$,
 $\l\in S_{1/2}[0,\pi]$.
In particular,
the estimate is valid
 for $\l\in S_{1/2}[0,\delta]$, $z\in\Gamma_\l^-$.
    \item $\nu=*$ and
    $z\in D_*^\delta(\l)$,
 $\l\in S_{1/2}[\delta,\pi]$.
In particular,
the estimate is valid
 for
$\l\in S_{1/2}[\frac{\pi}{2}-\frac{\delta}{2},\pi]$, $z\in\Gamma_\l$.
\end{description}
\end{theorem}

{\bf Proof of Theorem \ref{Basic:Svoistva}.}
We present the proof
only for  $A_0$  when $\l\in S_{1/2}[0,\pi]$ and
$z\in S(-\pi,\frac{\pi}{3})\cap D_Z^\delta(\l)$;
for other cases the proof is analogous.
Note that existence of a solution of (\ref{EqinZ})
implies its  analyticity in $D_Z^0(\l)$
(by Fuchs theorem, see Ch.5 \S 3 in \cite{OlverBook}).

The scheme of the proof is as follows.
First we replace (\ref{EqinZ}) by an equivalent integral equation.
Since its integrand is analytic,
we choose the integration path to be a curve $\Upsilon_\varphi(z)$,
given by Definition~\ref{DefinitionOfYpsilon} (see Appendix~A).
Using the estimate of the integrand on these curves,
 standard iteration scheme yields the required estimate of
  $z$-derivative of the solution.

For any  $z\in S(-\pi,\frac{\pi}{3})\cap D_Z^\delta(\l)$ there exist
 a contour $\Upsilon_\varphi(z)$ with
 $|\varphi|\le\frac{\pi}{3}-\frac{\delta}{3}$,
which lies within the same domain:
$\Upsilon_\varphi(z)\subset S(-\pi,\frac{\pi}{3})\cap D_Z^\delta(\l)$.
Since $V_0$ is analytic
(see its definition (\ref{EffPotentials})),
$u$ solves (\ref{EqinZ}) if it is a solution of the
 integral equation
 \begin{equation}
\label{IntegralEqinZ}
  u(z,\l)=\Ai(z)+
\int_{\Upsilon_\varphi(z)}
   J_0(z,s)V_0(s,\l) u(s,\l)\,ds,\quad
  \quad z\in S(-\pi,\tfrac{\pi}{3})\cap D_Z^\delta(\l),
\end{equation}
where
 $
J_0(z,s)=\Ai(s)\Bi(z)-\Ai(z)\Bi(s)
$.
Here $\int_{\Upsilon_\varphi(z)} f(s) ds$ denotes the complex line integral
of $f$ along the infinite curve ${\Upsilon_\varphi(z)}$.
We have to treat (\ref{IntegralEqinZ}) as a formal equation;
it is justified below using standard iteration technique.
We rewrite the last equation in terms of
$v=e^{\frac{2}{3}z^{\frac{3}{2}}}u$:
 \begin{equation}
 \label{IEpodkrychennoe}
  v(z,\l)=\ai(z)+
  \int_{\Upsilon_\varphi(z)}
   J(z,s)V_0(s,\l)v(s,\l)\,ds,\quad  \ai(z)\equiv
   \Ai(z)e^{\frac{2}{3}z^{\frac{3}{2}}},\quad
 \end{equation}
 where
 $
J(z,s)=J_0(z,s)
e^{\frac{2}{3}(z^{\frac{3}{2}}-s^{\frac{3}{2}})}
 $.
To provide
continuity as $\arg\l\downarrow0$, everywhere below
we require that $z^{\frac{3}{2}}$ takes its
values on the lower side of the cut $(-\infty,0]$
for $z<0$.  (Since  for
$\Im\l>0$   the curve $\Gamma_\l$ lies in the lower half-plane $\Im z<0$.)
 By Definition~\ref{DefASmall} and (\ref{BfromA}),
we have
\begin{equation}
\label{PodkrutkaU}
J(z,s)
=
-2\pi ie^{-i\frac{\pi}{3}}\left(
\ai(z)\ai(s\omega)-
e^{\frac{4}{3}(z^{\frac{3}{2}}-s^{\frac{3}{2}})} \ai(z\omega)\ai(s)\right),
\quad
\omega=e^{\frac{2\pi i}{3}}.
\end{equation}
Set $v_0(z)=a(z)$ and consider the iterations
$$
v_n(w,\l)=
\int_{\Upsilon_\varphi(w)}
J(w,s)V_0(s,\l)v_{n-1}(s,\l)ds
\quad\text{for }\quad
w\in \Upsilon_\varphi(z) .
$$
We estimate $v_n$ in terms of the majorizing functions
$
b_n(w,\l)=\sup\limits_{s\in\Upsilon_\varphi(w)}(1+|s|)^{\frac{1}{4}}|v_n(s,\l)|
$,
defined on
$\Upsilon_\varphi(z)$.
By \refLem{Upsilon:Decreasing}{PropertiesOfUpsilon},
$b_n(w,\l)$ is non-decreasing on $\Upsilon_\varphi(z)$.
By (\ref{AsAiry}),
 \begin{equation}\label{aiest}
|\ai(z)|\le \frac{\co}{(1+|z|)^{\frac{1}{4}}}, \quad 
|\ai'(z)|\le \frac{\co}{(1+|z|)^{\frac{5}{4}}}, \qquad |\arg z|\le
\pi-\ve,\quad \forall\ve>0. 
\end{equation}
Therefore, using (\ref{EstV0-z}) and \refLem{Upsilon:Monotonicity}{PropertiesOfUpsilon},
we obtain
$$
|J(w,s)V_0(s,\l)|
\le
\frac{\co}{(1+|w|)^{\frac{1}{4}}(|\l|^{\frac{4}{3}}+|s|^2)(1+|s|)^{\frac{1}{4}}}
\quad
\text{for}
\quad
s\in\Upsilon_\varphi(w)
$$
uniformly in $w\in S(-\pi,\frac{\pi}{3})\cap D_Z^\delta(\l)$.
Thus
$$
b_{n+1}(w)
\le
\co
\int_{\Upsilon_\varphi(w)}
\frac{b_n(s) |ds|}{(|\l|^{\frac{4}{3}}+|s|^2)(1+|s|)^{\frac{1}{2}}}
\le\frac{\co}{|\l|^{\frac{2}{3}}}
\int_{\Upsilon_\varphi(z)}
\frac{b_n(s) |ds|}{(1+|s|)^{\frac{3}{2}}},
\qquad w\in\Upsilon_\varphi(z),
$$
where we denote by $\int_G f(s) \,|ds|$  the line integral of $f$
along $G$ with respect to the arc length
$|ds|=\sqrt{(dx)^2+(dy)^2}$.
By (\ref{aiest}), we have $b_0\le\co$, so the integrals converge absolutely.
Using the induction principle, we  obtain
\begin{equation}
\label{EstVSolution}
(1+|w|)^{\frac{1}{4}}|v_n(w,\l)|
\le
b_n(w,\l)
\le
\frac{1}{n!}
\left(
\frac{\co}{|\l|^{\frac{2}{3}}}
\int_{\Upsilon_\varphi(w)}
\frac{ |ds|}{(1+|s|)^{\frac{3}{2}}}
\right)^n,
\qquad
w\in\Upsilon_\varphi(z).
\end{equation}
By \refLem{Upsilon:PowIntegration}{PropertiesOfUpsilon},
the integral in parenthesis is bounded. Hence,
the  series $\sum_{n=0}^\infty v_n$ converges uniformly and absolutely
and its sum $v=\sum_{n=0}^\infty v_n$ solves (\ref{IEpodkrychennoe}).
Thus $u=e^{-\frac{2}{3}z^{\frac{3}{2}}}v$ solves (\ref{EqinZ}).
Moreover, (\ref{EstVSolution}) implies that for
$\l\in S_{1/2}[0,\pi]$
and
$z\in S(-\pi,\frac{\pi}{3})\cap D_Z^\delta(\l)$
\begin{equation}\label{Estv}
|v(z,\l)|
\le
\frac{1}{(1+|z|)^{\frac{1}{4}}}
\sum\limits_{n=0}^\infty
b_n(z,\l)
\le
\frac{\co}{(1+|z|)^{\frac{1}{4}}}
.
\end{equation}
Now  show that $u$ has the asymptotics (\ref{DefA_0}).
 We write $v$  as
\begin{equation}\label{DecomposeSolutionV}
v(z,\l)=
a(z)+a(z)I_p-a(z\omega)e^{\frac{4}{3}z^{\frac{3}{2}}}I_e,
\end{equation}
where
$$
I_p=-2\pi ie^{-i\frac{\pi}{3}}
\int\limits_{\Upsilon_\varphi(z)}V_0(s,\l)a(s\omega)v(s,\l)ds,
\quad
I_e=2\pi ie^{-i\frac{\pi}{3}}
\int\limits_{\Upsilon_\varphi(z)}
e^{-\frac{4}{3}s^{\frac{3}{2}}}
V_0(s,\l)a(s)v(s,\l)ds.
$$
By (\ref{EstV0-z}), (\ref{aiest}),  (\ref{Estv}) and
\refLem{Upsilon:ExpIntegration}{PropertiesOfUpsilon},
\begin{equation}\label{EstI_e}
 |I_e|
 \le
 \frac{\co}{|\l|^{\frac{2}{3}}}
 \frac{|e^{-\frac{4}{3}z^{\frac{3}{2}}}|}{(1+|z|)^2}
 \qquad
 \text{for}
 \quad
 z\in S[-\pi+\tfrac{\delta}{3},\tfrac{\pi}{3}]\cap D_Z^\delta(\l).
\end{equation}
In order to estimate $I_p$ we observe that it is given by an integral
with an analytic integrand. So for
$z\in S(-\pi+\frac{4\vt}{3}+\ve, \frac{\pi}{3})$ we deform the
integration path to the ray
$\{s:\arg s=\arg z, |s|\ge|z|\}$, which lies within
$S(-\pi,\frac{\pi}{3})\cap D_Z^\delta(\l)$. Thus using
 (\ref{aiest}), (\ref{EstV0-z}) and (\ref{Estv}),
we obtain
\begin{equation}\label{EstI_p}
|I_p|
\le
\frac{\co}{(1+|z|)^{\frac{3}{2}}},
\quad
|I_p|
\le
\frac{\co}{|\l|^{\frac{2}{3}}(1+|z|)^{\frac{1}{2}}}
\qquad
 \text{for}
 \quad
z\in S[-\pi+\tfrac{4\vt}{3}+\ve,\tfrac{\pi}{3}].
\end{equation}
Using $\ai(z)=   \Ai(z)e^{\frac{2}{3}z^{\frac{3}{2}}}$,
 (\ref{DecomposeSolutionV}), (\ref{EstI_e}) and  (\ref{EstI_p}),
we obtain for
$u=e^{-\frac{2}{3}z^{\frac{3}{2}}}v$
the estimates
$$
|u(z,\l)-\Ai(z)|
\le
\co
\frac{|e^{-\frac{2}{3}z^{\frac{3}{2}}}|}{(1+|z|)^{\frac{3}{2}+\frac{1}{4}}},
\qquad
|u(z,\l)-\Ai(z)|
\le
\co
\frac{|e^{-\frac{2}{3}z^{\frac{3}{2}}}|}{|\l|^{\frac{2}{3}}(1+|z|)^{\frac{1}{2}+\frac{1}{4}}}
$$
uniformly in $\l\in S_{1/2}[0,\pi]$ and $z\in
S[-\pi+\frac{4\vt}{3}+\ve, \frac{\pi}{3}]$.
By the first estimate and (\ref{AsAiry}), $u$ satisfies
the first asymptotics in  (\ref{DefA_0}).
By the second estimate, the error term in (\ref{DefA_0})
can be replaced by $O(|\l|^{-\frac{2}{3}}|z|^{-\frac{1}{2}})$.

Taking $z$-derivative of (\ref{DecomposeSolutionV}) and using
(\ref{EstI_e}) and (\ref{EstI_p}), we obtain for
$\l\in S_{1/2}[0,\pi]$ and $z\in
S[-\pi+\frac{4\vt}{3}+\ve, \frac{\pi}{3}]$ the estimates
$$
|\partial_zu(z,\l)
-\Ai'(z)|
\le
\co
\frac{|e^{-\frac{2}{3}z^{\frac{3}{2}}}|}{(1+|z|)^{\frac{3}{2}-\frac{1}{4}}},
\quad
|\partial_zu(z,\l)
-\Ai'(z)|
\le
\co
\frac{|e^{-\frac{2}{3}z^{\frac{3}{2}}}|}{|\l|^{\frac{2}{3}}(1+|z|)^{\frac{1}{2}-\frac{1}{4}}},
$$
so that $u$
satisfies the second asymptotics in
 (\ref{DefA_0}) and its error term
can be replaced by $O(|\l|^{-\frac{2}{3}}|z|^{-\frac{1}{2}})$.

We demonstrated that $u$ has the asymptotics (\ref{DefA_0}), hence
 $A_0=u=e^{-\frac{2}{3}z^{\frac{3}{2}}}v$ and  $a_0=v$.
 Now prove (\ref{MainEstimateForASmall}) for $a_0$.
By (\ref{aiest}), (\ref{EstV0-z}), (\ref{Estv}) and
\refLem{Upsilon:PowIntegration}{PropertiesOfUpsilon},
we conclude that $I_p$ is uniformly bounded for
$\l\in S_{1/2}[0,\pi]$ and
$z\in S[-\pi+\frac{\delta}{3},\frac{\pi}{3}]\cap D_Z^\delta(\l)$.
Taking $z$-derivative of (\ref{DecomposeSolutionV}) and estimating
$|\partial_z v(z,\l)- a'(z)|$
using  (\ref{aiest}), (\ref{EstI_e}) and boundedness of $I_p$, we
obtain (\ref{MainEstimateForASmall}) for $\nu=0$.
 $\blacksquare$

Next we relate the solutions $A_\nu$ to $\psi(\pm x,\l)$ and
$\psi(\pm ix,-\l)$. The next Corollary follows from
 the asymptotics (\ref{DefA_0}--\ref{DefA_*}) for $A_\nu$,
 (\ref{DefinitionParabolic CylFunct})
 for $\psi$ and
$z_\l(x)= (\frac{3}{4})^{\frac{2}{3}} x^{\frac{4}{3}}
(1+O(x^{-1}))$ for $|\arg x-\vt|<\pi$ (see (\ref{DefXi}--\ref{DefZofX})).
Note that multiplication by $\phi$ (or $\phi^{-1}$),
which has different values on the upper and the lower
sides of the cut $\mathbb{R}_-$, annihilates with the similar
behavior of $A_\nu/\sqrt{z_\l'}$.
\begin{corollary}
  Let $|\l|\ge\frac{1}{2}$,
  $\l\in\mathbb{C}\setminus\mathbb{R}_-$ and $z_\l(x)\in D_Z(\l)$.
 Then
    \begin{align}
    \label{AoToPsi}
A_0)& &
        \psi(x,\l)
        &=\phi(\l)\frac{A_0(z_\l(x),\l)}{\sqrt{z_\l'(x)}}
        &
         \text{for} \quad
         \l&\in S_{1/2}(-\pi,\pi),
         \\
A_-) & &
        \psi(-ix,-\l)
        &=
        \frac{2^{\frac{3}{2}}\pi}{\phi(\l)}
        e^{i\frac{\pi}{12}}e^{i\frac{\pi}{4}\l}
        \frac{A_-(z_\l(x),\l)}{\sqrt{z_\l'(x)}}
&
         \text{for} \quad
         \l&\in S_{1/2}[-\pi+\delta,\pi],
\\
A_+) & &
        \psi(ix,-\l)
        &=
        \frac{2^{\frac{3}{2}}\pi}{\phi(\l)}
        e^{-i\frac{\pi}{12}}e^{-i\frac{\pi}{4}\l}
        \frac{A_+(z_\l(x),\l)}{\sqrt{z_\l'(x)}}
        &
         \text{for} \quad
\l&\in S_{1/2}[-\pi,\pi-\delta],
\\
\label{A*ToPsi}
A_*)& &
        \psi(-x,\l)
        &=
        e^{\pm i\frac{\pi}{6}}e^{\mp i\frac{\pi}{2}\l}\phi(\l)
        \frac{A_*(z_\l(x),\l)}{\sqrt{z_\l'(x)}}
&
         \text{for} \quad
\pm\arg\l&\in [\delta,\pi], |\l|\ge1/2,
      \end{align}
      where
  $\phi(\l)$ is given by (\ref{Example}).
\end{corollary}

 {\bf Proof of Theorem~\ref{Basic:InVariableX}.}
We give the proof only for $\psi(x,\l)$ and $\Im\l\ge0$; for other cases
 the proof is analogous.
By (\ref{AoToPsi}), we have
$
\phi(\l)a_0(\zlx,\l)=
\psi(x,\l)\sqrt{\zlpx}
\exp({\frac{2}{3}\zlx^{\frac{3}{2}}})
$. Taking $x$-derivative of this identity and using (\ref{DefZofX}), we obtain
$$
\psi'(x,\l)
+\psi(x,\l)\sqrt{x^2-\l}
=
\phi(\l)
e^{-\l\xi(\tfrac{x}{\sqrt{\l}})}
\left(
\sqrt{z_\l'}\frac{\partial}{\partial z}a_0(z_\l,\l)
-\frac{z_\l''}{2(z_\l')^{\frac{3}{2}}}a_0(z_\l,\l)
\right),
$$
where we write $z_\l$ in place of $z_\l(x)$
for brevity.
Using the estimate (\ref{MainEstimateForASmall})
of Theorem~\ref{Basic:Svoistva}, we have
\begin{equation}\label{EstCombPsiPsiPrime-1}
|\psi'(x,\l)
+\psi(x,\l)\sqrt{x^2-\l}|
\le C_\delta |\phi(\l)e^{-\l\xi(\tfrac{x}{\sqrt{\l}})}|
\left(
\frac{\sqrt{|z_\l'|}}{1+|z_\l|^{\frac{5}{4}}}
+
\left|
\frac{z_\l''}{(z_\l')^{\frac{3}{2}}}
\right|
\frac{1}{1+|z_\l|^{\frac{1}{4}}}
\right)
\end{equation}
for $\zlx\in D_0^\delta(\l)\subset D_Z(\l)$. Using (\ref{DefXi}--\ref{DefZofX}),
we conclude that
for $x\in D_X(\l)$ (that is, $z_\l(x)\in D_Z(\l)$)
\begin{equation}\label{EstZlAndZlPrime}
|\zlx|\le \co \frac{|x^2-\l|}{|\l|^{\frac{1}{3}}+|x^2-\l|^{\frac{1}{3}}},
\quad
|z_\l'(x)|\le \co (|\l|^{\frac{1}{6}}+|x^2-\l|^{\frac{1}{6}}),
\quad
\left|
\frac{z_\l''}{z_\l'}
\right|
\le\co
\frac{|z_\l'|}{1+|z_\l|}.
\end{equation}
Now (\ref{EstCombPsiPsiPrime-1}) and (\ref{EstZlAndZlPrime})
yield (\ref{EstPsi-0}), as required.
$\blacksquare$

\begin{corollary}[Symmetries]
Let $z\in D_Z(\l)$. Then
\begin{equation}
\overline{A_0(\bar{z},\l)} = A_0(z,\bar{\l}),
\qquad
\overline{A_*(\bar{z},\l)} = A_*(z,\bar{\l}),
\qquad
\overline{A_{\pm}(\bar{z},\l)} = A_{\mp}(z,\bar{\l}).
\end{equation}
\end{corollary}

Next we use (\ref{AoToPsi}--\ref{A*ToPsi})
and the identity
$
\phi(-\l)\phi(\l)=2^{\frac{3}{2}}\pi e^{\pm i\frac{\pi}{4}\l}$
for
$
\pm\Im\l>0
$ to rewrite the
connection formulas (\ref{ConnFormIX}--\ref{ConnFormMinusX})
 in terms of $A_\nu$.

\begin{corollary}[Connection formulas]
Let $A_\alpha$ denote $A_\alpha(z,\l)$ for $\alpha=0,+,-,*$, where
$z\in D_Z(\l)$. Then
\begin{align}
\label{ConnFormA0}
  A_0 &=
  \tfrac{2\sqrt{\pi}}{\phi^2(\l)}
  \Gamma(\tfrac{\l+1}{2})
  \left[
  e^{-i\frac{\pi}{3}}A_++e^{i\frac{\pi}{3}}A_-
  \right],
  &\text{for} \quad
  -\pi<\arg\lambda
  &<\pi,
   \\
  A_\pm&= \frac{e^{\mp i\frac{\pi}{2}\l}}{2\cos \frac{\pi}{2}\l}
  \frac{1}{\frac{2\sqrt{\pi}}{\phi^2(\l)}
  \Gamma({\tfrac{\l+1}{2}})}
  \left[
   e^{\pm i\pi\l} e^{\pm i\frac{\pi}{3}}A_0 +A_*
  \right],
  &\text{for} \quad
  0<\pm\arg\lambda
  &<\pi,
   \\
A_\mp&= \frac{e^{\mp i\frac{\pi}{2}\l}}{2\cos \frac{\pi}{2}\l}
  \frac{1}{
  \frac{2\sqrt{\pi}}{\phi^2(\l)}
  \Gamma({\tfrac{\l+1}{2}})
  }
  \left[
    e^{\pm i\frac{\pi}{3}}A_* +e^{\mp i\frac{\pi}{3}}A_0
  \right],
 &\text{for} \quad
  0<\pm\arg\lambda
  &<\pi,
  \\
\label{ConnFormA*}
A_* &=
  \tfrac{2\sqrt{\pi}}{\phi^2(\l)}
  \Gamma(\tfrac{\l+1}{2})
  \left[
    e^{\pm i\pi\l}e^{\mp i\frac{\pi}{3}}A_\mp +A_\pm
  \right],
  &\text{for} \quad
  0<\pm\arg\lambda
  &<\pi.
\end{align}
\end{corollary}

Next we find the Wronskians  $W\{f,g\}=fg'-f'g$ for the solutions $A_\nu$,
using the asymptotics (\ref{DefA_0}--\ref{DefA_*})
and the connection formulas (\ref{ConnFormA0}--\ref{ConnFormA*}).
\begin{corollary}
For $0<\arg\l<\pi$ we have
\begin{equation}\label{}
 W\{A_0,A_{\pm}\} =
 \frac{e^{\mp i\frac{\pi}{6}}}{2\pi},
 \qquad
W\{A_*,A_+\} =
-\frac{e^{ i\frac{\pi}{6}}}{2\pi}e^{i\pi\l},
\qquad
W\{A_*,A_-\} =
  \frac{e^{i\frac{\pi}{2}}}{2\pi},
\end{equation}
\begin{equation}\label{}
W\{A_0,A_*\} =
\frac{e^{- i\frac{\pi}{6}}}{\pi}
e^{i\frac{\pi}{2}\l}\cos(\tfrac{\pi\l}{2})
\tfrac{2\sqrt{\pi}}{\phi^2(\l)}
\Gamma(\tfrac{\l+1}{2}),
\qquad
W\{A_-,A_+\} =
\frac{e^{- i\frac{\pi}{2}}}{2\pi}
\frac{1}{
  \frac{2\sqrt{\pi}}{\phi^2(\l)}
  \Gamma({\tfrac{\l+1}{2}})
  }.
\end{equation}
\end{corollary}

\section{Appendix A}\label{AppendixYpsilon}
\renewcommand{\theequation}{A.\arabic{equation}}
\setcounter{equation}{0}
We define the family of curves $\Upsilon_\varphi(z)$
and study its properties. We use it
in the proof of
Theorem~\ref{Basic:Svoistva}.

\begin{definition}
\label{DefinitionOfYpsilon}
For a complex point  $ z\in S(-\pi,\pi)$ and an angle
$\varphi\in[-\frac{\pi}{3},\frac{\pi}{3}]$, satisfying
$|\arg z-\varphi|\le\frac{2\pi}{3}$, we set
\begin{equation}\label{DefUpsilon}
\Upsilon_\varphi(z)
=
\left\{
s\in S|[\arg z,\varphi]| :
\Im(se^{-i\varphi})^{\frac{3}{2}}=\Im(ze^{-i\varphi})^{\frac{3}{2}},
\Re(se^{-i\varphi})^{\frac{3}{2}}\ge\Re(ze^{-i\varphi})^{\frac{3}{2}}
\right\},
\end{equation}
 where
$S|[\arg z,\varphi]|$ denotes the sector $S[\arg z,\varphi]$ if
$\arg z\le\varphi$ and $S[\varphi,\arg z]$ otherwise.
\end{definition}
  The curve $\Upsilon_\varphi(z)$ is asymptotic to the ray $\arg
s=\varphi$.
If $\arg z=\varphi$, then $\Upsilon_\varphi(z)$ degenerates into
the ray $e^{i\varphi}[|z|,\infty)$. If
$|\arg z-\varphi|=\frac{2\pi}{3}$, then $\Upsilon_\varphi(z)$ degenerates into
the sum of the interval
$e^{i\arg z}[0,|z|]$ and the ray
$e^{i\varphi}\mathbb{R}_+$.
The  curves $\arg(ze^{i\varphi})^{3/2}=$const are schematically presented on
Fig.\ref{Ypsilon} a).

\begin{lemma}
\label{PropertiesOfUpsilon}
\begin{enumerate}
    \item \label{Upsilon:Decreasing}
    If $w\in \Upsilon_\varphi(z)$, then
$\Upsilon_\varphi(w)\subset \Upsilon_\varphi(z)$.
    \item \label{Upsilon:Monotonicity}
    If $s\in\Upsilon_\varphi(z)\setminus\{z\}$ and
    $|\varphi|<\frac{\pi}{3}$, then
$\Re(s^{\frac{3}{2}}-z^{\frac{3}{2}})=
\cos\frac{3\varphi}{2}\cdot
\Re((se^{-i\varphi})^{\frac{3}{2}}-
(ze^{-i\varphi})^{\frac{3}{2}})
>0
$.
    \item \label{Upsilon:ExpIntegration}
    Let $\alpha\in\mathbb{R}$ and $\delta>0$. Then for
    $|\arg z|\le\pi-\frac{2\delta}{3}$ there exists $\varphi$ such that
     $|\varphi|\le\frac{\pi}{3}-\frac{\delta}{3}$,
$|\arg z-\varphi|\le\frac{2\pi}{3}-\frac{\delta}{3}$ and
\begin{equation}\label{ExpAlongYpsilon}
\int_{\Upsilon_\varphi(z)}
\frac{|e^{-\frac{4}{3}{s}^{\frac{3}{2}}}|}{(1+| s|)^{\alpha}}\,|ds|
\le C_\delta
\frac{|e^{-\frac{4}{3}z^{\frac{3}{2}}}|}{(1+|z|)^{\alpha+\frac{1}{2}}},
\qquad
\end{equation}
where $C_\delta$ is independent of $z$ and $\varphi$.
\item \label{Upsilon:PowIntegration}
Let $\alpha>1$. Then  the integral
$
\int_{\Upsilon_\varphi(z)}
\frac{|ds|}{(1+| s|)^{\alpha}}
$
is bounded uniformly in $z$.
\end{enumerate}
\end{lemma}
{\bf Proof.}
\ref{Upsilon:Decreasing} is  evident. To prove \ref{Upsilon:Monotonicity}
observe that
$s^{\frac{3}{2}}=e^{\frac{3i\varphi}{2}}(se^{-i\varphi})^{\frac{3}{2}}$,
where we take the principal value of non-integer powers on
$\mathbb{C}\setminus\mathbb{R}_-$.
Hence,
$\Re s^{\frac{3}{2}}
=
\cos\tfrac{3\varphi}{2}\cdot\Re(se^{-i\varphi})^{\frac{3}{2}}
+
\sin\tfrac{3\varphi}{2}\cdot\Im(se^{-i\varphi})^{\frac{3}{2}}$,
where the last term  is constant for all
$s\in\Upsilon_\varphi(z)$. Subtracting the same formula with $s=z$
yields the result.

Now prove \ref{Upsilon:ExpIntegration}.
Let $x=\Re (ze^{-i\varphi})^{\frac{3}{2}}$,
$y=\Im (ze^{-i\varphi})^{\frac{3}{2}}$.
Parametrize
the curve $\Upsilon_\varphi(z)$ by
$s(t)=e^{i\varphi}(t+iy)^{\frac{2}{3}}$, $t\in[x,\infty)$.
By \ref{Upsilon:Monotonicity},
$
\Re s^{\frac{3}{2}}
=
t\cos\frac{3\varphi}{2}
+
y\sin\frac{3\varphi}{2}
$;
using
$|ds|=\frac{2}{3}\frac{dt}{|t+iy|^{\frac{1}{3}}}$,
we obtain
\begin{equation}\label{IntAlongR+(w)}
 \int_{\Upsilon_\varphi(z)}
\frac{|e^{-\frac{4}{3}{s}^{\frac{3}{2}}}|}{(1+| s|)^{\alpha}}\,|ds|
\le C e^{-y\frac{4}{3}\sin\frac{3\varphi}{2}}I,
\qquad
I=
\int_x^\infty
\frac{e^{-x\frac{4}{3}\cos\frac{3\varphi}{2}} dt}{(1+|t+iy|)^{\frac{2\alpha}{3}}|t+iy|^{\frac{1}{3}}},
\end{equation}
where we used
$(1+|s|)^{-\alpha}\le\co(1+|t+iy|)^{-\frac{2\alpha}{3}}$.
It remains to prove that
\begin{equation}\label{EstI}
  I\le
  \frac{\co \cdot e^{-\varepsilon x}}{(1+|w|)^{\frac{2}{3}(\alpha+\frac{1}{2})}},
  \quad\hbox{where}\quad
  w=x+iy,
   |\arg w|\le\pi-\frac{\delta}{2}
  \quad\hbox{and}\quad
\varepsilon=\tfrac{4}{3}\cos\tfrac{3\varphi}{2}.
\end{equation}
For $|w|\le1$ this is evident.
For $|w|>1$ and $|\arg w|\le\frac{\pi}{2}$ we use
$|t+iy|\ge |w|$  to obtain
$
I\le
\frac{C}{(1+|w|)^{\frac{2}{3}(\alpha+\frac{1}{2})}}
\int_x^\infty e^{-\varepsilon t} dt
$. Since $\ve\ge\frac{4}{3}\sin\frac{\delta}{2}$,
this yields (\ref{EstI}). It remains to consider $|w|>1$ and
$\frac{\pi}{2}\le|\arg w|\le\pi-\frac{\delta}{2}$. Using
 $\sin\frac{\delta}{2}\le|w|\sin\frac{\delta}{2}\le |y|$, we have
$$
I\le
\frac{1}{|y|^{\frac{2}{3}(\alpha+\frac{1}{2})}}
\int_x^\infty e^{-\varepsilon t} dt
\le
\frac{3}{4
(\sin\frac{\delta}{2})^{\frac{2}{3}(\alpha+\frac{1}{2})+1}
}
\frac{
e^{-\varepsilon x}
}{
|w|^{\frac{2}{3}(\alpha+\frac{1}{2})}}
\le
  \co \frac{e^{-\varepsilon x}}{(1+|w|)^{\frac{2}{3}(\alpha+\frac{1}{2})}},
$$
which also implies (\ref{EstI}). Substituting (\ref{EstI}) into
(\ref{IntAlongR+(w)})
completes the proof.

To prove \ref{Upsilon:PowIntegration} we use
the same  parametrization $s(t)$ as in the proof of
\ref{Upsilon:ExpIntegration}. This gives
\begin{equation}\label{IntAlongR+Pow}
 \int_{\Upsilon_\varphi(z)}
\frac{|ds|}{(1+| s|)^{\alpha}}\,|ds|
\le
\int_x^\infty
\frac{\tfrac{2}{3}dt}{(1+|t+iy|^{\frac{2}{3}})^\alpha|t+iy|^{\frac{1}{3}}}
\le
\frac{2}{\alpha-1}.
\blacksquare
\end{equation}

\section{Appendix B: Integration along $\Gamma_\l$}
\renewcommand{\theequation}{B.\arabic{equation}}
\setcounter{equation}{0}

In this section we estimate the integrals of $\frac{1}{(1+|z|)^\alpha}$
and
$\frac{|e^{\pm z^{\frac{3}{2}}}|}{(1+|z|)^{\alpha}}$.  We show that for $\delta\le|\arg \l|\le\pi$
the integrals along the family of
curves $\Gamma_\l(z)$ allow the same estimates as the integrals along
$\mathbb{R}_+$. For $|\arg \l|\le\delta$ this is true only for
$z\in\Gamma_\l^+$.
The  integrals along
$\Gamma_\l^-$ for $|\arg \l|\le\delta$ allow the same estimate as those along
$[-|\l|^{\frac{2}{3}},0]$.

 For any continuous function
$f$  on a smooth curve $G$
we denote the usual complex line integral by $\int_G f(s)\,ds$.
 We denote by $\int_G f(s) \,|ds|$  the line integral of $f$
along $G$ with respect to the arc length
$|ds|=\sqrt{(dx)^2+(dy)^2}$. For integration along the infinite curve
$G$
we use the same notation
 $\int_{G} f(s)\,ds=\lim\limits_{R\to\infty} \int\limits_{G\cap\{s:|s|\le R\}} f(s)\,ds$
 for absolutely converging integrals.

Now we formulate the main result of this section.

\begin{theorem}
\label{EstExpInt}
Let $|\l|\ge1/2$ and $\alpha\in\mathbb{R}$. Fix $\delta\in(0,\frac{\pi}{5})$
and assume
that either a)
$z\in\Gamma_\l$,
$ \delta\le|\arg\l|\le\pi$ or b) $z\in\Gamma_\l^+$. Then  the following
estimates are fulfilled:
\begin{align}
\label{Exp-2} \int_{\Gamma_\l(z)}
\frac{|e^{-\frac{4}{3}{s}^{\frac{3}{2}}}|}{(1+| s|)^{\alpha}}\,|ds|
&\le C
\frac{|e^{-\frac{4}{3}z^{\frac{3}{2}}}|}{(1+|z|)^{\alpha+\frac{1}{2}}},
\\
\label{ExpGrow}
\int_{\Gamma_\l(w,z)}
\frac{|e^{\frac{4}{3}{s}^{\frac{3}{2}}}|}{(1+| s|)^{\alpha}}\,|ds|
&\le C
\frac{|e^{\frac{4}{3}z^{\frac{3}{2}}}|}{(1+|z|)^{\alpha+\frac{1}{2}}},
\quad\hbox{where}\quad
\left\{\begin{array}{lcl}
  w=z_0  & for & \delta\le|\arg\l|\le\pi, \\
  w=z_*  & for & |\arg\l|\le\delta ,\\
\end{array}
\right.
\\
\label{pow-2}
\int_{\Gamma_\l(z)}
\frac{|ds|}{(1+|s|)^{\alpha}}
&\le \frac{C}{(1+|z|)^{\alpha-1}},
\qquad
\alpha>1,
\\
\label{DDecayPowerG+}
\int_{\Gamma_\l(z)}
\frac{(1+|s|)^{-1}|ds|}{(1+|s||\l|^{-\frac{2}{3}})^{\alpha}}
&\le
C\frac{\alpha^{-1}+\ln(1+2|\l|)}{(1+|z||\l|^{-\frac{2}{3}})^{\alpha}},
\qquad \alpha>0,
\end{align}
where  $C$ is independent of $\l$ and $z$.
\end{theorem}

\begin{theorem}
\label{PowLem} Let $|\l|\ge1/2$.  Then
 the following
estimates are fulfilled:
\begin{equation}
\label{pow-3Le1}
\int_{\Gamma_\l^-}
\frac{|ds|}{(1+|s|)^{\alpha}}\, \le
\left\{
\begin{array}{ll}
  {C}(1-\alpha)^{-1}|\l|^{\frac{2}{3}(1-\alpha)}\quad &
\hbox{\rm for}\quad 0\le\alpha<1,\\[1ex]
  \co \log (1+2|\l|) \quad &
  \hbox{\rm for}\quad \alpha=1,\\[1ex]
{C}{(\alpha-1)^{-1}} \quad &
 \hbox{\rm for}\quad \alpha>1,\\
\end{array}
\right.
\end{equation}
\begin{equation}\label{EIV0}
\int_{\Gamma_\l} \frac{|ds|}{|\l|^{\frac{4}{3}}+|s|^2}
\le\frac{C}{|\l|^{\frac{2}{3}}},
\end{equation}
where  $C$ is independent of $\l$ and $z$.
\end{theorem}

We consider only the case $\Im\l\ge0$; for $\Im\l\le0$ the proof is analogous.
As a prerequisite for the proof we estimate
$\int_{\Gamma_\l(w,z)}|f(s)||ds|$
for the  cases
a) and b) of the hypothesis of  Theorem~\ref{EstExpInt}.
Similarly we estimate
$\int_{\Gamma_\l^-}|f(s)||ds|$
for $0\le\arg\l\le\delta$.

\begin{figure}
\centering
\includegraphics[width=165mm]{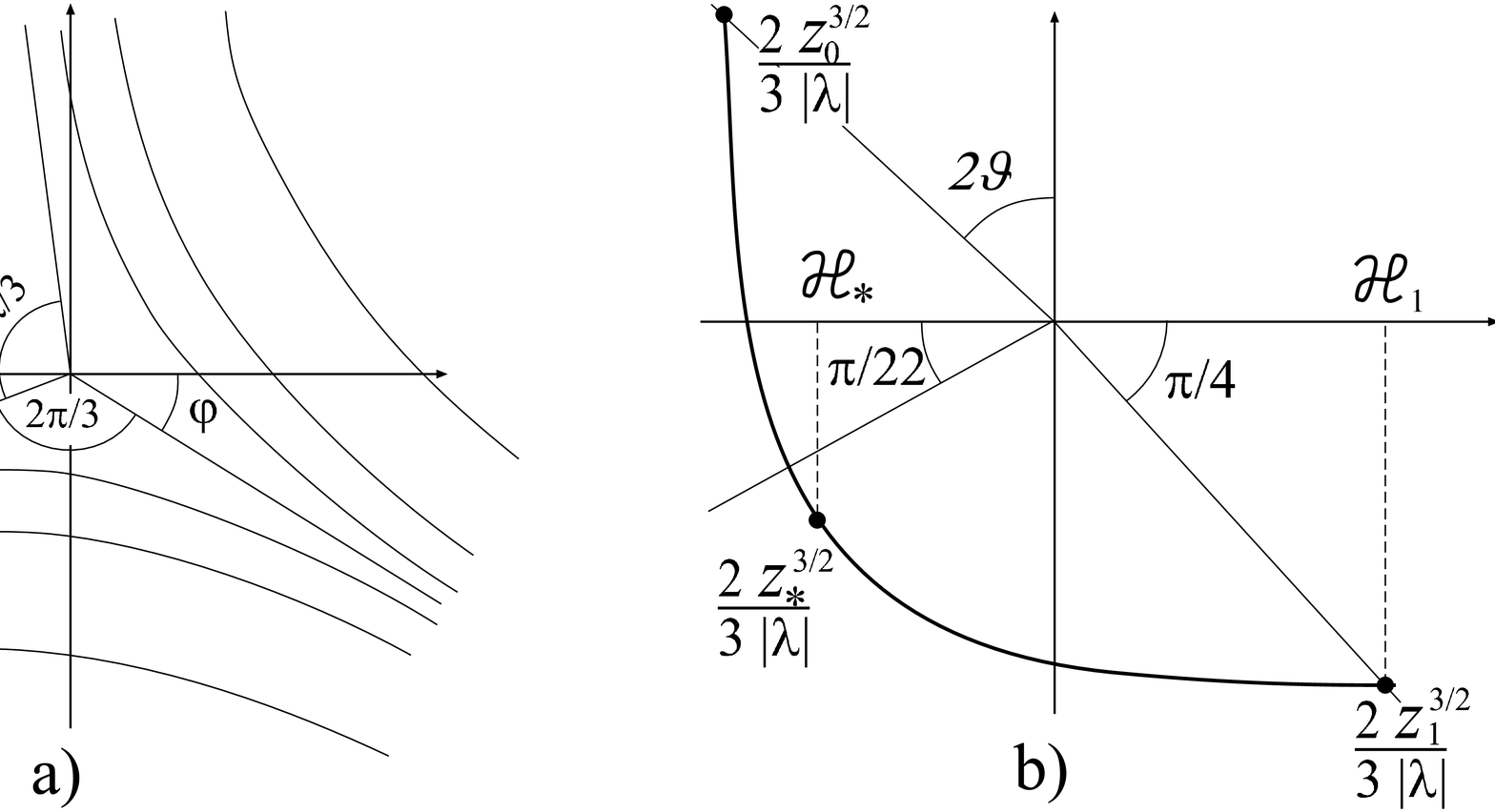}\\
\caption[short caption here]{
a) the curves $\arg(ze^{i\varphi})^{\frac{3}{2}}=$const
in the sector
$S[-\varphi-\frac{2\pi}{3},-\varphi+\frac{2\pi}{3}]$; \\
b) the curve $e^{2i\vt}\xi(e^{-i\vt}\mathbb{R}_+)$.  }
\label{Ypsilon}  
\end{figure}

Introduce a convenient parametrization of
$\Gamma_\l$ for $\delta<\arg\l\le\pi$ and of $\Gamma_\l^+$ for
$0\le\arg\l\le\delta$.
By definition of $\Gamma_\l=z_\l(\mathbb{R}_+)$,
the mapping $z_\l(\cdot)$  already gives the parametrization by $x\in\mathbb{R}_+$.
Define the new parameter
$\varkappa$ as a function of $x$ by $\varkappa=\Re( e^{2i\vt}\xi(\frac{x}{\sqrt{\l}}))$.
This is a smooth one-to-one mapping, since
by \refLem{A7.2-2}{23}, $\Re( e^{2i\vt}\xi(\frac{x}{\sqrt{\l}}))$ is strictly
increasing in $x$. We also introduce
 the function $v_\vt$, which maps
the real part of $e^{2i\vt}\xi(re^{-i\vt})$ ($r\ge0$) to its imaginary part:
for $\varkappa=\Re( e^{2i\vt}\xi(re^{-i\vt}))$ we set
 $v_\vt(\varkappa)=\Im( e^{2i\vt}\xi(re^{-i\vt}))$.
The curve
$e^{2i\vt}\xi(e^{-i\vt}\mathbb{R}_+)$ is
 schematically presented on Fig.\ref{Ypsilon} b).
The parametrization of $\Gamma_\l$ in terms of $\varkappa$ is given by
\begin{equation}\label{SOfKappa}
    s_\l(\varkappa)=|\l|^{\frac{2}{3}}
    \left[
\tfrac{3}{2}(\varkappa+iv_\vt(\varkappa))
    \right]^{\frac{2}{3}},
    \quad
\text{so that}
\quad
|ds|=\sqrt{1+\left(\frac{dv_\vt(\varkappa)}{d\varkappa} \right)^2}
\frac{|\l|d\varkappa}{\sqrt{|s_\l(\varkappa)|}}.
\end{equation}
By (\ref{DefEta}), (\ref{DefZofX}) and (\ref{SOfKappa}),
for a point $z_\l(x)=s_\l(\vk)$ on $\Gamma_\l$ we have
\begin{equation}\label{TwoParametrizations}
\tfrac{2}{3}\frac{z_\l^{\frac{3}{2}}(x)}{|\l|}
=
e^{2i\vt}\xi(t)=\varkappa+i v_\vt(\varkappa)
=
\tfrac{2}{3}\frac{s_\l^{\frac{3}{2}}(\vk)}{|\l|},
\qquad
\text{where}
\quad
t=\frac{x}{\sqrt{\l}}=re^{-i\vt}, \quad r\ge0.
\end{equation}
Let us show that
\begin{equation}\label{EstimateDvDkappa}
 \left|\frac{dv_\vt(\varkappa)}{d\varkappa}\right|\le C
\qquad
\text{in cases}
\qquad
\begin{array}{lll}
\textrm{a)} & s_\l(\varkappa)\in\Gamma_\l, &  \delta<\arg\l\le\pi,  \\[1ex]
\textrm{b)} & s_\l(\varkappa)\in\Gamma_\l^+,  & 0\le\arg\l\le\delta. \\
\end{array}
\end{equation}
For $\varkappa=\Re( e^{2i\vt}\xi(t))$, where $t=re^{-i\vt}$ and  $r\ge0$,
we have
\begin{equation}\label{ExplicitDvDKappa}
\frac{dv_\vt(\varkappa)}{d\varkappa}= \frac{ {\partial_r}
\Im
\left(e^{2i\vt}\xi(re^{-i\vt}) \right)}{ {\partial_r} \Re
\left(
e^{2i\vt}\xi(re^{-i\vt})\right)
 }= \tan \arg \left(e^{2i\vt}{\partial_r}
\xi(t)\right)
=\tan \arg \left(
e^{i\vt}\sqrt{t^2-1}
\right).
\end{equation}

In case a)
 by \refLem{A7.2-1}{Geom-0}, we have $\arg{\partial_r}\xi(t)
\in [-\frac{\pi}{2}-\vt,-2\vt)$.
Hence
$\arg \left(e^{2i\vt}{\partial_r} \xi(t)\right)
\in
[-\frac{\pi}{2}+\frac{\delta}{2}, 0)$ and
$
\left|\frac{d v_\vt(\varkappa)}{d\varkappa}\right|\le\cot\frac{\delta}{2}
$
uniformly for $\delta\le\arg \l\le \pi$.

In case b)  by definition (\ref{DefDelta}), we have
 $0\le\vt<\frac{\pi}{10}$, so that
(\ref{Argxi(z*)}) (equivalent to \refLem{A7.3-1}{DivG2})
  implies $-\pi-\vt\le\arg\xi(t)$.
 Therefore by \refLem{A7-5}{Geom-0},
$\arg\left( e^{2i\vt}{\partial_r}
\xi(t)\right)\in[-\frac{\pi}{3}+\frac{\vt}{6},\frac{\vt}{2}]$.
Hence
$\left|
\frac{d v_\vt(\varkappa)}{d\varkappa}
\right|
\le\tan\frac{\pi}{3}
$,
as required.

Thus for
$\Gamma_\l(z_1,z_2)\subset\Gamma_\l$, $\delta<\arg\l\le\pi$
and for $\Gamma_\l(z_1,z_2)\subset\Gamma_\l^+$,
$0\le\arg\l\le\delta$ the estimate (\ref{EstimateDvDkappa}) implies
\begin{equation}\label{UToKappaEstimate}
\int\limits_{\Gamma_\l(z_1,z_2)}|f(s)||ds|\le
\co |\l|
\int\limits_{\varkappa_1}^{\varkappa_2}
\frac{|f(s_\l(\varkappa))|}{\sqrt{|s_\l(\varkappa)|}}\,d\varkappa,
\qquad
\text{where}
\quad
\varkappa_{1,2}=\Re \frac{2}{3}\frac{z_{1,2}^{\frac{3}{2}}}{|\l|}.
\end{equation}
By \refLem{A.7.2-3a}{23}, the last inequality in \refLem{A7.3-1}{DivG2} and (\ref{TwoParametrizations}),
we have
\begin{equation}\label{VTK}
-\pi+\frac{\pi}{22}
\le
\arg (\vk+iv_\vt(\vk))
\qquad\hbox{for}\quad
 s_\l(\vk)\in\Gamma_\l^+,
 \quad
  0<\arg\l\le\pi.
\end{equation}
By (\ref{ExplicitDvDKappa}), (\ref{TwoParametrizations}), \refLem{A7-6}{Geom-0} and
(\ref{DefDelta}),
we have
\begin{equation}\label{DerVTForPosKappa}
\left|\frac{dv_\vt(\vk)}{d\vk}\right|\le\tan(\tfrac{\pi}{4}-\tfrac{\pi}{15})<1
\qquad\text{if}
\quad
\vk\ge0,
\quad
0<\arg\l\le\delta.
\end{equation}

Now we estimate $\int_{\Gamma_\l^-}|f(s)||ds|$
for $0\le\arg\l\le\delta$.
We introduce the parametrization
of $\Gamma_\l^-$, symmetric in a sense to (\ref{SOfKappa}):
 now the imaginary part of $\frac{z^{\frac{3}{2}}}{|\l|}$ becomes the parameter.
 Let $x$ parametrizes $\Gamma_\l$ by $z=z_\l(x)$;
for $0\le\arg\l\le\delta$ define the new parameter
$\chi$ as a function of $x$  by $\chi(x)=\Im( e^{2i\vt}\xi(\frac{x}{\sqrt{\l}}))$.
This is a smooth one-to-one mapping, since
by  \refLem{A7.2-2}{23}, $\Im( e^{2i\vt}\xi(\frac{x}{\sqrt{\l}}))$
is strictly decreasing.
We introduce the function $u_\vt$,
 which maps
the imaginary part of $e^{2i\vt}\xi(re^{-i\vt})$ ($r\ge0$) to its real part:
for $\chi=\Im( e^{2i\vt}\xi(re^{-i\vt}))$ we set
$u_\vt(\chi)=\Re( e^{2i\vt}\xi(re^{-i\vt}))$.
The parametrization of $\Gamma_\l^-$ in terms of $\chi$ is
\begin{equation}\label{SOfChi}
    w_\l(\chi)=|\l|^{\frac{2}{3}}
    \left[
\tfrac{3}{2}(u_\vt(\chi)+i\chi)
    \right]^{\frac{2}{3}},
    \quad
|dw|=\sqrt{1+\left(\frac{du_\vt(\chi)}{d\chi} \right)^2}
\frac{|\l|\,d\chi}{\sqrt{|w_\l(\chi)|}}.
\end{equation}
Let us show that $\frac{du_\vt}{d\chi}$ is uniformly bounded.
For $\chi=\Im( e^{2i\vt}\xi(t))$, where $t=re^{-i\vt}$ and  $r\ge0$,
we have
$$
\frac{du_\vt(\chi)}{d\chi}= \frac{ {\partial_r}
\Re
\left(e^{2i\vt}\xi(re^{-i\vt}) \right)}{ {\partial_r} \Im
\left(
e^{2i\vt}\xi(re^{-i\vt})\right)
 }= \cot \arg\left( e^{2i\vt}{\partial_r}
\xi(t)\right).
$$
By \refLem{A7.3-3}{DivG2} and
(\ref{DefDelta}),
$\arg \left(e^{2i\vt}{\partial_r}
\xi(t)\right)\in [-\frac{\pi}{2},-\frac{\pi}{6}]$.
Therefore $\left|\frac{du_\vt(\chi)}{d\chi}\right|$
is uniformly bounded
and for $0\le\arg\l\le\delta$ we have
\begin{equation}\label{UToChiEstimate}
\int\limits_{\Gamma_\l^-}|f(s)||ds|\le
\co |\l|
\int\limits_{\chi_*}^{\chi_0}
\frac{|f(w_\l(\chi))|}{\sqrt{|w_\l(\chi)|}}\, d\chi,
\qquad
\text{where}
\quad
\chi_{0,*}=\Im \frac{2}{3}\frac{z_{0,*}^{\frac{3}{2}}}{|\l|}
\end{equation}
and
$\chi_0=\frac{\pi}{4}\cos 2\vt$.

\begin{lemma}\label{Z1}
Let $\l\in S_{1/2}[0,\pi]$. For each $\l$ define
 $\vk_0$, $\vk_*$ and $\vk_1>0$ by $s_\l(\vk_{0,*})=z_{0,*}$
and $v_\vt(\vk_1)=-\vk_1$.
 Then $\vk_*\in[\vk_0,\vk_1]$ and
there exist a positive number $C$, independent of $\l$ and $z$, such that
\begin{enumerate}
    \item \label{Z1-1}
    $|s_\l(\vk)|\le C|z_*|$ for $\delta\le\arg\l\le\pi$ and
    $\vk\in[\vk_0,\vk_1]$,
    \item \label{Z1-2}
    $|s_\l(\vk)|^{\frac{3}{2}}\le C|\Im z_*^{\frac{3}{2}}|$ for
    $0\le\arg\l\le\delta$ and
    $\vk\in[\vk_*,\vk_1]$,
    \item \label{Z1-3}
    $|s_\l(\vk)|^{\frac{3}{2}}\le \frac{3}{\sqrt{2}}|\l|
    \vk$ for $0\le\arg\l\le\pi$ and
    $\vk\in[\vk_1,\infty)$.
\end{enumerate}
\end{lemma}
\noindent {\bf Proof.}
Our main instrument is the relation (\ref{TwoParametrizations}).
Together with (\ref{DerVTForPosKappa}) it implies uniqueness of $\vk_1$.
Using  \refLem{A.7.2-3a}{23} and \refLem{A7.3-1}{DivG2},
we conclude that $\vk_*\in[\vk_0,\vk_1]$.
By  (\ref{TwoParametrizations})
 and \refLem{A7.2-2}{23}, $v_\vt(\varkappa)$ is non-increasing.
 The  points $\vk_{0,*,1}$
are schematically presented on Fig.\ref{Ypsilon} b).

\ref{Z1-1}.
Write $t$ such that $\l^{\frac{2}{3}}\eta(t)=s_\l(\vk)$
in the form (\ref{AngularXi}):
$t=re^{-i\vartheta}=1+\eta e^{-i\varphi}$, $\varphi\in[0,\pi]$,
$\eta\ge0$.
By \refLem{A.7.2-3a}{23} and the definition of $z_1$,
 we have $\arg\xi(t)\le-\frac{\pi}{4}-2\vt$. Now using \refLem{A7-2}{Geom-0}
 we obtain $-\frac{\pi}{2}\le\vt-\varphi\le-\frac{\pi}{6}$, so that the identity
 $
r=\frac{\sin\varphi}{\sin(\varphi-\vt)}
$ yields $|t|\le\frac{1}{\sin\frac{\pi}{6}}$.
Applying (\ref{TwoParametrizations}),
we conclude that
$
\frac{2}{3}s_\l(\vk)^{\frac{3}{2}}
\subset
\l\xi(\{
t:|t|\sin\frac{\pi}{6}\le1
\})
$
for $\vk\in[\vk_0,\vk_1]$, hence
$|s_\l(\vk)|\le\co|\l|^{\frac{2}{3}}$.
\refLem{A7.5-2}{SvaGForAll} gives
 $|\l|^{\frac{2}{3}}\sin\frac{\delta}{2}\le|z_*|$.
This yields $|s_\l(\vk)|\le\co|z_*|$, as required.

\ref{Z1-2}. It suffice to consider $\arg\l\neq0$.
The estimates
(\ref{EstimateDvDkappa}), (\ref{VTK}) implies
$|v_\vt(0)|\le \co|v_\vt(\vk_*)|$ and
(\ref{DerVTForPosKappa}) implies
 $|v_\vt(\vk_1)|\le \co|v_\vt(0)|$,
 so that $|v_\vt(\vk_1)|\le \co|v_\vt(\vk_*)|$.
 By (\ref{TwoParametrizations}), (\ref{DefDelta}), \refLem{A7.3-1}{DivG2} in the form
(\ref{Argxi(z*)}) and definition of $\vk_1$,
 we have
$
 -\frac{3\pi}{4}-\frac{\pi}{8}\le\arg(\vk+iv_\vt(\vk))\le-\frac{\pi}{4}
\qquad\text{for}\quad
\vk\in[\vk_*,\vk_1]
$.
Taking into account non-increasing of $v_\vt(\vk)$, we conclude that
$|\vk|\le|v_\vt(\vk)|\cot\frac{\pi}{8}$ and therefore
$|\vk+iv_\vt(\vk)|\le\co|v_\vt(\vk_*)|$ for $\vk\in[\vk_*,\vk_1]$.
By (\ref{TwoParametrizations}),
this is equivalent to
$|s_\l(\vk)|^{\frac{3}{2}}\le C|\Im z_*^{\frac{3}{2}}|$, as required.

\ref{Z1-3}. By definition of $\vk_1$, (\ref{TwoParametrizations}),
  and \refLem{A.7.2-3a}{23},
 we have $s_\l(\vk)^{\frac{3}{2}}\in S[-\frac{\pi}{4},0]$
 for $\vk\in[\vk_1,\infty)$,
 as required.
 $\blacksquare$
Note that
$\vk_0=\Re\frac{2}{3}\frac{z_0^{\frac{3}{2}}}{|\l|}=-\frac{\pi}{4}\sin2\vt$.

\noindent {\bf Proof of Theorem~\ref{EstExpInt}.}
First we prove (\ref{Exp-2}). Consider the case $0\le\arg\l\le\pi$,
$z\in\Gamma_\l^+$. By \refLem{A7.3-2}{DivG2},
for $s_\l(\varkappa)\in\Gamma_\l^+$ the function
$|s_\l(\varkappa)|$ is non-decreasing. Thus using (\ref{SOfKappa}),
(\ref{UToKappaEstimate}),
$|\l\varkappa|^{\frac{2}{3}}\le\sqrt{|s_\l(\varkappa)|}$
for $|z|\le1$ and $\sqrt{|z|}\le\sqrt{|s_\l(\varkappa)|}$ for $|z|>1$,
we obtain
$$
\int\limits_{\Gamma_\l(z)}
\frac{|e^{-\frac{4}{3}{s}^{\frac{3}{2}}}|}{(1+| s|)^{\alpha}}\,|ds|
\le
\frac{C|\l|}{(1+|z|)^{\alpha}}
\int\limits_{\Re \frac{2}{3}\frac{z^{\frac{3}{2}}}{|\l|}}^\infty
\frac{e^{-2|\l|\varkappa }\, d\varkappa}{\sqrt{|s_\l(\varkappa)|}}
\le C
\frac{|e^{-\frac{4}{3}z^{\frac{3}{2}}}|}{(1+|z|)^{\alpha+\frac{1}{2}}}.
$$

Consider the case $\delta\le\arg\l\le\pi$,
$z\in\Gamma_\l^-$. By \ref{A7.5-2} and \ref{A7.5-3} of Lemma~\ref{SvaGForAll},
$|z_*|$ is bounded away from zero and $|z|\le C|z_*|$. So using
(\ref{SOfKappa}),
(\ref{UToKappaEstimate}) and the definition  of $z_*$ (\ref{Z-star}) we obtain
$$
\int\limits_{\Gamma_\l(z)}
\frac{|e^{-\frac{4}{3}{s}^{\frac{3}{2}}}|}{(1+| s|)^{\alpha}}\,|ds|
\le \frac{C|\l|}{(1+|z_*|)^{\alpha+\frac{1}{2}}}
\int\limits_{\Re \frac{2}{3}\frac{z^{\frac{3}{2}}}{|\l|}}^\infty
e^{-2|\l|\varkappa }d\varkappa
\le C
\frac{|e^{-\frac{4}{3}z^{\frac{3}{2}}}|}{(1+|z|)^{\alpha+\frac{1}{2}}}
$$
as required. This completes the proof of (\ref{Exp-2}).

Consider (\ref{ExpGrow}). Introduce the notations
$a=\Re\frac{2}{3}\frac{z^{\frac{3}{2}}}{|\l|}$,
$b=\Re\frac{2}{3}\frac{w^{\frac{3}{2}}}{|\l|}$ (here and below
we omit dependence of $\l$).
Applying (\ref{UToKappaEstimate})
we obtain
\begin{equation}\label{DefIForPosExp}
\int\limits_{\Gamma_\l(w,z)}
\frac{|e^{\frac{4}{3}{s}^{\frac{3}{2}}}|}{(1+| s|)^{\alpha}}\,|ds|
\le C I(b,a),
\quad\text{where}\quad
I(b,a)=\int\limits_{b}^{a}
\frac{e^{2|\l|\vk}|\l| d\vk}{(1+|s_\l(\vk)|)^\alpha\sqrt{|s_\l(\vk)|}}.
\end{equation}
Let $\delta<\arg\l\le\pi$. Suppose that $z\in\Gamma_\l(z_0,z_1)$, where
$z_0$ is given by (\ref{Defz0zE}) and $z_1=s_\l(\vk_1)$ for $\vk_1$
is defined in Lemma~\ref{Z1}.
By \ref{A7.5-2} and \ref{A7.5-3} of Lemma~\ref{SvaGForAll},
$|z_*|$ is bounded away from zero and $|z|\le C|z_*|$. Therefore,
$$
I(b,a)\le
I(\vk_0,a)
\le
\frac{1}{(1+|z_*|)^\alpha|z_*|^{\frac{1}{2}}}
\int_{\vk_0}^{a}e^{2|\l|\vk}|\l|d\vk
\le
\co \frac{e^{2|\l|a}}{(1+|z_*|)^{\alpha+\frac{1}{2}}}
\le
C
\frac{|e^{\frac{4}{3}z^{\frac{3}{2}}}|}{(1+|z|)^{\alpha+\frac{1}{2}}}.
$$
 By (\ref{DefIForPosExp}), the last estimate
proves (\ref{ExpGrow}) for $z\in\Gamma_\l(z_0,z_1)$.

Now suppose that $z\in\Gamma_\l(z_1,\infty)$. By
\refLem{Z1-3}{Z1}, $a>0$ and $|z|\le\co(|\l|a)^{\frac{2}{3}}$.
Therefore,
$$
I(\vk_0,a/2)\le
\int_{\vk_0}^{a/2}e^{2|\l|\vk}|\l|d\vk
\le
e^{|\l|a}
\le
\co \frac{e^{2|\l|a}}{(1+(|\l|a)^{\frac{2}{3}})^{\alpha+\frac{1}{2}}}
\le
C\frac{|e^{\frac{4}{3}z^{\frac{3}{2}}}|}{(1+|z|)^{\alpha+\frac{1}{2}}},
$$
where we used $e^{-t}\le\frac{\co}{(1+t)^\alpha}$ for $t\ge0$.
Using \refLem{A7.5-2}{SvaGForAll}, \refLem{Z1-3}{Z1} and $a>0$,
we obtain
$$
I(a/2,a)
\le
C \frac{e^{2|\l|a}}{(|\l|a)^{\frac{2}{3}(\alpha+\frac{1}{2})}}
\le
C\frac{|e^{\frac{4}{3}z^{\frac{3}{2}}}|}{(1+|z|)^{\alpha+\frac{1}{2}}}.
$$
By (\ref{DefIForPosExp}), the two last displayed formulas prove (\ref{ExpGrow})
for $z\in\Gamma_\l(z_1,\infty)$.

Now let $0<\arg\l\le\delta$ (for $\l>0$ the proof is by direct calculation). For $|z|\le1$ the
result is evident, so we consider the case $|z|>1$.
Suppose that $z\in\Gamma_\l(z_*,z_1)$.
By \refLem{A7.5-2}{SvaGForAll}
and \refLem{Z1-2}{Z1}, $0<|v_\vt(\vk_*)\l|^{\frac{2}{3}}\le|z_*|\le|z|\le\co|v_\vt(\vk_*)\l|^{\frac{2}{3}}$.
Hence,
$$
I(\vk_*,a)
\le
\co \int_{\vk_*}^{a}
\frac{e^{2|\l|\vk}|\l|\,d\vk}{|v_\vt(\vk_*)\l|^{\frac{2}{3}\alpha+\frac{1}{3}}}
\le
C\frac{|e^{\frac{4}{3}z^{\frac{3}{2}}}|}{|z|^{\alpha+\frac{1}{2}}}
\le
C\frac{|e^{\frac{4}{3}z^{\frac{3}{2}}}|}{(1+|z|)^{\alpha+\frac{1}{2}}}.
$$
By (\ref{DefIForPosExp}), this proves (\ref{ExpGrow})
for $z\in\Gamma_\l(z_*,z_1)$.

It remains to consider $z\in\Gamma_\l(z_1,\infty)$.
By definition of $z_1$, $a\ge0$.
Using the substitution $t=|\vk\l|^{\frac{2}{3}}$
and \refLem{Z1-3}{Z1}, we
 have
$$
I(\vk_*,a)
\le
\int\limits_{-\infty}^{a}
\frac{\co\cdot e^{2|\l|\vk}|\l|d\vk}{(1+|\l\vk|^{\frac{2}{3}})^\alpha|\l\vk|^{\frac{1}{3}}}
\le
\int\limits_{-\infty}^{(a|\l|)^{\frac{2}{3}}}
\frac{\co\cdot e^{2t^{\frac{3}{2}}}}{(1+|t|)^\alpha}dt
\le
\frac{\co\cdot e^{2|\l|a}}{(1+(a|\l|)^{\frac{2}{3}})^{\alpha+\frac{1}{2}}}
\le\co
\frac{|e^{\frac{4}{3}z^{\frac{3}{2}}}|}{(1+|z|)^{\alpha+\frac{1}{2}}}.
$$
By (\ref{DefIForPosExp}), this proves (\ref{ExpGrow}).

Next we  prove (\ref{pow-2}). Consider the case
$0\le\arg\l\le\pi$,
$z\in\Gamma_\l^+$. By \refLem{A7.2-2}{23}, $v_\vt(\varkappa)$ is
non-increasing. By (\ref{TwoParametrizations}) and
the last estimate in
\refLem{A7.3-1}{DivG2},
$v_\vt(\varkappa)\le0$. Thus using
(\ref{SOfKappa}),
 (\ref{UToKappaEstimate}) and the substitution
$t=|\varkappa|^{\frac{2}{3}}$ we obtain
\begin{equation}\label{LocI}
 \int\limits_{\Gamma_\l(z)}
\frac{|ds|}{(1+|s|)^{\alpha}} \le
\frac{\co\cdot I}{|\l|^{\frac{2}{3}(\alpha-1)}},\quad
I=
\int\limits_{a}^\infty
\frac{\frac{2}{3}|\varkappa|^{-\frac{1}{3}} d\varkappa}{
(|\varkappa|^{\frac{2}{3}}+|v|^{\frac{2}{3}}+\varepsilon)^\alpha}
=\int\limits_{|a|^{\frac{2}{3}}\sign a}^\infty
\frac{dt}{(|t|+|v|^{\frac{2}{3}}+\varepsilon)^\alpha},
\end{equation}
where $a=\Re \frac{2}{3}\frac{z^{\frac{3}{2}}}{|\l|}$,
$v=\Im \frac{2}{3}\frac{z^{\frac{3}{2}}}{|\l|}$,
$\ve=(\frac{3}{2}|\l|)^{-\frac{2}{3}}$.
For $a\ge0$ direct calculation yields
$I\le\frac{\co}{\alpha-1}(|a|^{\frac{2}{3}}+|v|^{\frac{2}{3}}+\varepsilon)^{1-\alpha}$;
by (\ref{LocI}), this gives (\ref{pow-2}). For $a<0$ we have
\begin{equation}
\label{PowEstG+}
I\le
2\int_0^\infty
\frac{dt}{(t+|v|^{\frac{2}{3}}+\varepsilon)^\alpha}
\le
\frac{2/(\alpha-1)}{(|v|^{\frac{2}{3}}+\varepsilon)^{\alpha-1}}
\le
\frac{\co/(\alpha-1)}{(|a|^{\frac{2}{3}}+|v|^{\frac{2}{3}}+\varepsilon)^{\alpha-1}},
\end{equation}
where the last estimate follows from (\ref{VTK}).
Now (\ref{PowEstG+}) and (\ref{LocI}) again give (\ref{pow-2}).

It remains to consider the case
$\delta\le\arg\l\le\pi$,
$z\in\Gamma_\l^-$.
By (\ref{TwoParametrizations}) and \refLem{A7.5-2}{SvaGForAll},
we have
$(\frac{2}{3})^{2}(\sin\frac{\delta}{2})^{3}
\le
\varkappa^2+v_\vt^2(\varkappa)$. Thus using
(\ref{SOfKappa}),
(\ref{UToKappaEstimate})  and \refLem{A7.5-3}{SvaGForAll}, we have
\begin{equation}\label{PowEstBigArgLambda}
\int_{\Gamma_\l(z)}
\frac{|ds|}{(1+|s|)^{\alpha}}
\le
\frac{\co}{|\l|^{\frac{2}{3}(\alpha-1)}}
\int\limits_{-\infty}^\infty
\frac{|\varkappa|^{-\frac{1}{3}} d\varkappa}{(\varkappa^2+v_\vt^2(\varkappa))^{\frac{\alpha}{3}}}
\le
\frac{\co}{|\l|^{\frac{2}{3}(\alpha-1)}}
\le
\frac{\co}{(1+|z|)^{\alpha-1}},
\end{equation}
as required.

Now prove (\ref{DDecayPowerG+}). The proof is  similar to that of
(\ref{pow-2}). Consider the case $0\le\arg\l\le\pi$,
$z\in\Gamma_\l^+$. By \refLem{A7.2-2}{23}, $v_\vt(\varkappa)$ is
non-increasing; by (\ref{TwoParametrizations}) and the last estimate in
\refLem{A7.3-1}{DivG2},
  we have $v_\vt(\varkappa)\le0$. Thus using (\ref{SOfKappa}),
 (\ref{UToKappaEstimate}) and the substitution
$t=|\varkappa|^{\frac{2}{3}}$ we obtain
$$
\int\limits_{\Gamma_\l(z)}
\frac{|ds|}{(1+|s|)(1+|s||\l|^{-\frac{2}{3}})^{\alpha}}
\le
\co J,
\quad
J=
\int\limits_{|a|^{\frac{2}{3}}\sign a}^\infty
\frac{dt}{(\ve+|v|^{\frac{2}{3}}+|t|)(1+|v|^{\frac{2}{3}}+|t|)^\alpha},
$$
where $a=\Re \frac{2}{3}\frac{z^{\frac{3}{2}}}{|\l|}$,
$v=\Im \frac{2}{3}\frac{z^{\frac{3}{2}}}{|\l|}$,
$\ve=(\frac{3}{2}|\l|)^{-\frac{2}{3}}$.
For $a\ge0$ direct estimate of the last integral  gives (\ref{DDecayPowerG+}).
For $a<0$ we expand the integration range to
$(-\infty,\infty)$ and use symmetry of the integrand.
This gives
$$
J\le2
\int\limits_{|v|^{\frac{2}{3}}}^\infty
\frac{dt}{(\ve+t)(1+t)^\alpha}
\le
2\frac{\alpha^{-1}+\ln(1+\ve^{-1})}{(1+|v|^{\frac{2}{3}})^\alpha}.
$$
Now we use (\ref{VTK}) to deduce
(\ref{DDecayPowerG+}) from the last inequality.

It remains to estimate the integral over
$\Gamma_\l(z,z_*)$ for
 $z\in\Gamma_\l^-$,
$\delta\le\arg\l\le\pi$.
By (\ref{TwoParametrizations})
and  \refLem{A7.5-3}{SvaGForAll},
 $\varkappa_0$ and $z\l^{-\frac{2}{3}}$
are bounded.
 So we use (\ref{SOfKappa}),
 (\ref{UToKappaEstimate}) and the substitution
$t=|\varkappa|^{\frac{2}{3}}$ to obtain
$$
\int_{\Gamma_\l(z)}
\frac{(1+|s|)^{-1}|ds|}{(1+|s||\l|^{-\frac{2}{3}})^{\alpha}}
\le \co
\int_0^{|\varkappa_0|^{\frac{2}{3}}}
\frac{dt}{\ve+t}
\le\co\ln(1+2|\l|)\le
C\frac{\alpha^{-1}+\ln(1+2|\l|)}{(1+|z||\l|^{-\frac{2}{3}})^{\alpha}}.
$$
 Combining
the last estimate with the result
for $z\in\Gamma_\l^+$ completes the proof.
$\blacksquare$

\textbf{Proof of Lemma~\ref{PowLem}.}
It is sufficient to consider $\l\in S_{1/2}[0,\pi]$.
First prove (\ref{pow-3Le1}); by (\ref{UToChiEstimate}), we have
\begin{equation}
\label{PowLem-c}
\int_{\Gamma_\l^-}
\frac{|ds|}{(1+ |s|)^\alpha}
\le
\co
\frac{I}{|\l|^{\frac{2}{3}(\alpha-1)}},
\qquad
I=\frac{3}{2}
\int_{\chi_*}^{\chi_0}
\frac{d\chi}{|\chi|^{\frac{1}{3}}(|\chi|^{\frac{2}{3}}+
 \varepsilon)^\alpha},
\end{equation}
where
$
\ve=(\frac{3}{2}|\l|)^{-\frac{2}{3}}$,
$\chi_0=\Im(e^{2i\vt}\xi(0))=\frac{\pi}{4}\cos2\vt$,
$\chi_*=\Im(e^{2i\vt}\xi(t_*))$.
By \refLem{A7.3-4}{DivG2}, $|\chi_*|$ is uniformly bounded;
by the last inequality in \refLem{A7.3-1}{DivG2}, $\chi_*\le0$.
The change of variable
$t=|\chi|^{\frac{2}{3}}\sign\chi$  in $I$ and further  direct estimate give (\ref{pow-3Le1}).

Now  prove (\ref{EIV0}). For $\delta\le\arg\l\le\pi$ we use
the parametrization of $\Gamma_\l$, given by (\ref{SOfKappa})
and the estimate (\ref{UToKappaEstimate}).
Similarly for $0\le\arg\l\le\delta$ we use (\ref{SOfKappa})
and  (\ref{UToKappaEstimate})
on $\Gamma_\l^+$ and
(\ref{SOfChi}), (\ref{UToChiEstimate}) on $\Gamma_\l^-$.
In the both cases we obtain
\begin{equation}
\label{PowLem-c1}
\int_{\Gamma_\l}
\frac{|ds|}{|\l|^{\frac{4}{3}}+|s|^2}
\le
\frac{\co}{|\l|^{\frac{2}{3}}}
\int_{-\infty}^\infty
\frac{dt}{|t|^{\frac{1}{3}}(1+|t|^{\frac{2}{3}})^{2}}
\le
\frac{\co}{|\l|^{\frac{2}{3}}}.
\qquad
\blacksquare
\end{equation}

\section{Acknowledgements}

The author is thankful to E.Korotyaev,
 A.Badanin and D.Chelkak for useful discussions.

\end{document}